\catcode`\@=11
\font\tensmc=cmcsc10      
\def\smc{\tensmc}
\def\pagewidth#1{\hsize= #1 }
\def\pageheight#1{\vsize= #1 }
\def\hcorrection#1{\advance\hoffset by #1 }
\def\vcorrection#1{\advance\voffset by #1 }
\def\wlog#1{}
\newif\iftitle@
\outer\def\title{\title@true\vglue 24\p@ plus 12\p@ minus 12\p@
   \bgroup\let\\=\cr\tabskip\centering
   \halign to \hsize\bgroup\tenbf\hfill\ignorespaces##\unskip\hfill\cr}
\def\endtitle{\cr\egroup\egroup\vglue 18\p@ plus 12\p@ minus 6\p@}
\outer\def\author{\iftitle@\vglue -18\p@ plus -12\p@ minus -6\p@\fi\vglue
    12\p@ plus 6\p@ minus 3\p@\bgroup\let\\=\cr\tabskip\centering
    \halign to \hsize\bgroup\smc\hfill\ignorespaces##\unskip\hfill\cr}
\def\endauthor{\cr\egroup\egroup\vglue 18\p@ plus 12\p@ minus 6\p@}
\outer\def\heading{\bigbreak\bgroup\let\\=\cr\tabskip\centering
    \halign to \hsize\bgroup\smc\hfill\ignorespaces##\unskip\hfill\cr}
\def\endheading{\cr\egroup\egroup\nobreak\medskip}
\outer\def\subheading#1{\medbreak\noindent{\tenbf\ignorespaces
      #1\unskip.\enspace}\ignorespaces}
\outer\def\proclaim#1{\medbreak\noindent\smc\ignorespaces
    #1\unskip.\enspace\sl\ignorespaces}
\outer\def\endproclaim{\par\ifdim\lastskip<\medskipamount\removelastskip
  \penalty 55 \fi\medskip\rm}
\outer\def\demo#1{\par\ifdim\lastskip<\smallskipamount\removelastskip
    \smallskip\fi\noindent{\smc\ignorespaces#1\unskip:\enspace}\rm
      \ignorespaces}
\outer\def\enddemo{\par\smallskip}
\newcount\footmarkcount@
\footmarkcount@=1
\def\makefootnote@#1#2{\insert\footins{\interlinepenalty=100
  \splittopskip=\ht\strutbox \splitmaxdepth=\dp\strutbox 
  \floatingpenalty=\@MM
  \leftskip=\z@\rightskip=\z@\spaceskip=\z@\xspaceskip=\z@
  \noindent{#1}\footstrut\rm\ignorespaces #2\strut}}
\def\footnote{\let\@sf=\empty\ifhmode\edef\@sf{\spacefactor
   =\the\spacefactor}\/\fi\futurelet\next\footnote@}
\def\footnote@{\ifx"\next\let\next\footnote@@\else
    \let\next\footnote@@@\fi\next}
\def\footnote@@"#1"#2{#1\@sf\relax\makefootnote@{#1}{#2}}
\def\footnote@@@#1{$^{\number\footmarkcount@}$\makefootnote@
   {$^{\number\footmarkcount@}$}{#1}\global\advance\footmarkcount@ by 1 }

\hyphenation{man-u-script man-u-scripts ap-pen-dix ap-pen-di-ces}
\hyphenation{data-base data-bases}
\ifx\amstexloaded@\relax\catcode`\@=13 
  \endinput\else\let\amstexloaded@=\relax\fi
\newlinechar=`\^^J
\def\eat@#1{}
\def\Space@.{\futurelet\Space@\relax}
\Space@. %
\newhelp\athelp@
{Only certain combinations beginning with @ make sense to me.^^J
Perhaps you wanted \string\@\space for a printed @?^^J
I've ignored the character or group after @.}
\def\futureletnextat@{\futurelet\next\at@}
{\catcode`\@=\active
\lccode`\Z=`\@ \lowercase
{\gdef@{\expandafter\csname futureletnextatZ\endcsname}
\expandafter\gdef\csname atZ\endcsname
   {\ifcat\noexpand\next a\def\next{\csname atZZ\endcsname}\else
   \ifcat\noexpand\next0\def\next{\csname atZZ\endcsname}\else
    \def\next{\csname atZZZ\endcsname}\fi\fi\next}
\expandafter\gdef\csname atZZ\endcsname#1{\expandafter
   \ifx\csname #1Zat\endcsname\relax\def\next
     {\errhelp\expandafter=\csname athelpZ\endcsname
      \errmessage{Invalid use of \string@}}\else
       \def\next{\csname #1Zat\endcsname}\fi\next}
\expandafter\gdef\csname atZZZ\endcsname#1{\errhelp
    \expandafter=\csname athelpZ\endcsname
      \errmessage{Invalid use of \string@}}}}
\def\atdef@#1{\expandafter\def\csname #1@at\endcsname}
\newhelp\defahelp@{If you typed \string\define\space cs instead of
\string\define\string\cs\space^^J
I've substituted an inaccessible control sequence so that your^^J
definition will be completed without mixing me up too badly.^^J
If you typed \string\define{\string\cs} the inaccessible control sequence^^J
was defined to be \string\cs, and the rest of your^^J
definition appears as input.}
\newhelp\defbhelp@{I've ignored your definition, because it might^^J
conflict with other uses that are important to me.}
\def\define{\futurelet\next\define@}
\def\define@{\ifcat\noexpand\next\relax
  \def\next{\define@@}%
  \else\errhelp=\defahelp@
  \errmessage{\string\define\space must be followed by a control 
     sequence}\def\next{\def\garbage@}\fi\next}
\def\undefined@{}
\def\preloaded@{}    
\def\define@@#1{\ifx#1\relax\errhelp=\defbhelp@
   \errmessage{\string#1\space is already defined}\def\next{\def\garbage@}%
   \else\expandafter\ifx\csname\expandafter\eat@\string
         #1@\endcsname\undefined@\errhelp=\defbhelp@
   \errmessage{\string#1\space can't be defined}\def\next{\def\garbage@}%
   \else\expandafter\ifx\csname\expandafter\eat@\string#1\endcsname\relax
     \def\next{\def#1}\else\errhelp=\defbhelp@
     \errmessage{\string#1\space is already defined}\def\next{\def\garbage@}%
      \fi\fi\fi\next}
\def\famzero{\fam\z@}

\def\arctan{\mathop{\famzero arctan}\nolimits}

\def\cos{\mathop{\famzero cos}\nolimits}

\def\det{\mathop{\famzero det}}

\def\exp{\mathop{\famzero exp}\nolimits}

\def\inf{\mathop{\famzero inf}}

\def\lim{\mathop{\famzero lim}}
\def\liminf{\mathop{\famzero lim\,inf}}
\def\limsup{\mathop{\famzero lim\,sup}}

\def\log{\mathop{\famzero log}\nolimits}
\def\max{\mathop{\famzero max}}
\def\min{\mathop{\famzero min}}

\def\sin{\mathop{\famzero sin}\nolimits}

\def\sup{\mathop{\famzero sup}}

\def\textfont@#1#2{\def#1{\relax\ifmmode
    \errmessage{Use \string#1\space only in text}\else#2\fi}}
\textfont@\rm\tenrm
\textfont@\it\tenit
\textfont@\sl\tensl
\textfont@\bf\tenbf
\textfont@\smc\tensmc
\let\ic@=\/
\def\/{\unskip\ic@}
\def\textfonti{\the\textfont1 }
\def\t#1#2{{\edef\next{\the\font}\textfonti\accent"7F \next#1#2}}
\let\B=\=
\let\D=\.
\def~{\unskip\nobreak\ \ignorespaces}
{\catcode`\@=\active
\gdef\@{\char'100 }}
\atdef@-{\leavevmode\futurelet\next\athyph@}
\def\athyph@{\ifx\next-\let\next=\athyph@@
  \else\let\next=\athyph@@@\fi\next}
\def\athyph@@@{\hbox{-}}
\def\athyph@@#1{\futurelet\next\athyph@@@@}
\def\athyph@@@@{\if\next-\def\next##1{\hbox{---}}\else
    \def\next{\hbox{--}}\fi\next}
\def\.{.\spacefactor=\@m}
\atdef@.{\null.}
\atdef@,{\null,}
\atdef@;{\null;}
\atdef@:{\null:}
\atdef@?{\null?}
\atdef@!{\null!}   
\def\srdr@{\thinspace}                     
\def\drsr@{\kern.02778em}
\def\sldl@{\kern.02778em}
\def\dlsl@{\thinspace}
\atdef@"{\unskip\futurelet\next\atqq@}
\def\atqq@{\ifx\next\Space@\def\next. {\atqq@@}\else
         \def\next.{\atqq@@}\fi\next.}
\def\atqq@@{\futurelet\next\atqq@@@}
\def\atqq@@@{\ifx\next`\def\next`{\atqql@}\else\def\next'{\atqqr@}\fi\next}
\def\atqql@{\futurelet\next\atqql@@}
\def\atqql@@{\ifx\next`\def\next`{\sldl@``}\else\def\next{\dlsl@`}\fi\next}
\def\atqqr@{\futurelet\next\atqqr@@}
\def\atqqr@@{\ifx\next'\def\next'{\srdr@''}\else\def\next{\drsr@'}\fi\next}
\def\flushpar{\par\noindent}
\def\textfontii{\the\textfont2 }
\def\{{\relax\ifmmode\lbrace\else
    {\textfontii f}\spacefactor=\@m\fi}
\def\}{\relax\ifmmode\rbrace\else
    \let\@sf=\empty\ifhmode\edef\@sf{\spacefactor=\the\spacefactor}\fi
      {\textfontii g}\@sf\relax\fi}   
\def\nonhmodeerr@#1{\errmessage
     {\string#1\space allowed only within text}}
\def\linebreak{\relax\ifhmode\unskip\break\else
    \nonhmodeerr@\linebreak\fi}
\def\allowlinebreak{\relax
   \ifhmode\allowbreak\else\nonhmodeerr@\allowlinebreak\fi}
\newskip\saveskip@
\def\nolinebreak{\relax\ifhmode\saveskip@=\lastskip\unskip
  \nobreak\ifdim\saveskip@>\z@\hskip\saveskip@\fi
   \else\nonhmodeerr@\nolinebreak\fi}
\def\newline{\relax\ifhmode\null\hfil\break
    \else\nonhmodeerr@\newline\fi}
\def\nonmathaerr@#1{\errmessage
     {\string#1\space is not allowed in display math mode}}
\def\nonmathberr@#1{\errmessage{\string#1\space is allowed only in math mode}}
\def\mathbreak{\relax\ifmmode\ifinner\break\else
   \nonmathaerr@\mathbreak\fi\else\nonmathberr@\mathbreak\fi}
\def\nomathbreak{\relax\ifmmode\ifinner\nobreak\else
    \nonmathaerr@\nomathbreak\fi\else\nonmathberr@\nomathbreak\fi}
\def\allowmathbreak{\relax\ifmmode\ifinner\allowbreak\else
     \nonmathaerr@\allowmathbreak\fi\else\nonmathberr@\allowmathbreak\fi}
\def\pagebreak{\relax\ifmmode
   \ifinner\errmessage{\string\pagebreak\space
     not allowed in non-display math mode}\else\postdisplaypenalty-\@M\fi
   \else\ifvmode\penalty-\@M\else\edef\spacefactor@
       {\spacefactor=\the\spacefactor}\vadjust{\penalty-\@M}\spacefactor@
        \relax\fi\fi}
\def\nopagebreak{\relax\ifmmode
     \ifinner\errmessage{\string\nopagebreak\space
    not allowed in non-display math mode}\else\postdisplaypenalty\@M\fi
    \else\ifvmode\nobreak\else\edef\spacefactor@
        {\spacefactor=\the\spacefactor}\vadjust{\penalty\@M}\spacefactor@
         \relax\fi\fi}
\def\newpage{\relax\ifvmode\vfill\penalty-\@M\else\nonvmodeerr@\newpage\fi}
\def\nonvmodeerr@#1{\errmessage
    {\string#1\space is allowed only between paragraphs}}
\def\smallpagebreak{\relax\ifvmode\smallbreak
      \else\nonvmodeerr@\smallpagebreak\fi}
\def\medpagebreak{\relax\ifvmode\medbreak
       \else\nonvmodeerr@\medpagebreak\fi}
\def\bigpagebreak{\relax\ifvmode\bigbreak
      \else\nonvmodeerr@\bigpagebreak\fi}
\newdimen\captionwidth@
\captionwidth@=\hsize
\advance\captionwidth@ by -1.5in
\def\caption#1{}
\def\topspace#1{\gdef\thespace@{#1}\ifvmode\def\next
    {\futurelet\next\topspace@}\else\def\next{\nonvmodeerr@\topspace}\fi\next}
\def\topspace@{\ifx\next\Space@\def\next. {\futurelet\next\topspace@@}\else
     \def\next.{\futurelet\next\topspace@@}\fi\next.}
\def\topspace@@{\ifx\next\caption\let\next\topspace@@@\else
    \let\next\topspace@@@@\fi\next}
 \def\topspace@@@@{\topinsert\vbox to 
       \thespace@{}\endinsert}
\def\topspace@@@\caption#1{\topinsert\vbox to
    \thespace@{}\nobreak
      \smallskip
    \setbox\z@=\hbox{\noindent\ignorespaces#1\unskip}%
   \ifdim\wd\z@>\captionwidth@
   \centerline{\vbox{\hsize=\captionwidth@\noindent\ignorespaces#1\unskip}}%
   \else\centerline{\box\z@}\fi\endinsert}
\def\midspace#1{\gdef\thespace@{#1}\ifvmode\def\next
    {\futurelet\next\midspace@}\else\def\next{\nonvmodeerr@\midspace}\fi\next}
\def\midspace@{\ifx\next\Space@\def\next. {\futurelet\next\midspace@@}\else
     \def\next.{\futurelet\next\midspace@@}\fi\next.}
\def\midspace@@{\ifx\next\caption\let\next\midspace@@@\else
    \let\next\midspace@@@@\fi\next}
 \def\midspace@@@@{\midinsert\vbox to 
       \thespace@{}\endinsert}
\def\midspace@@@\caption#1{\midinsert\vbox to
    \thespace@{}\nobreak
      \smallskip
      \setbox\z@=\hbox{\noindent\ignorespaces#1\unskip}%
      \ifdim\wd\z@>\captionwidth@
    \centerline{\vbox{\hsize=\captionwidth@\noindent\ignorespaces#1\unskip}}%
    \else\centerline{\box\z@}\fi\endinsert}
\mathchardef\prime@="0230
\def\prime{{{}\prime@{}}}
\def\prim@s{\prime@\futurelet\next\pr@m@s}

\def\,{\relax\ifmmode\mskip\thinmuskip\else\thinspace\fi}
\def\!{\relax\ifmmode\mskip-\thinmuskip\else\negthinspace\fi}
\def\frac#1#2{{#1\over#2}}

\def\binom#1#2{{#1\choose #2}}

\def\:{\nobreak\hskip.1111em{:}\hskip.3333em plus .0555em\relax}
\def\intic@{\mathchoice{\hskip5\p@}{\hskip4\p@}{\hskip4\p@}{\hskip4\p@}}
\def\negintic@
 {\mathchoice{\hskip-5\p@}{\hskip-4\p@}{\hskip-4\p@}{\hskip-4\p@}}
\def\intkern@{\mathchoice{\!\!\!}{\!\!}{\!\!}{\!\!}}
\def\intdots@{\mathchoice{\cdots}{{\cdotp}\mkern1.5mu
    {\cdotp}\mkern1.5mu{\cdotp}}{{\cdotp}\mkern1mu{\cdotp}\mkern1mu
      {\cdotp}}{{\cdotp}\mkern1mu{\cdotp}\mkern1mu{\cdotp}}}
\newcount\intno@             
\def\iint{\intno@=\tw@\futurelet\next\ints@} 
\def\iiint{\intno@=\thr@@\futurelet\next\ints@}
\def\iiiint{\intno@=4 \futurelet\next\ints@}
\def\idotsint{\intno@=\z@\futurelet\next\ints@}
\def\ints@{\findlimits@\ints@@}
\newif\iflimtoken@
\newif\iflimits@
\def\findlimits@{\limtoken@false\limits@false\ifx\next\limits
 \limtoken@true\limits@true\else\ifx\next\nolimits\limtoken@true\limits@false
    \fi\fi}
\def\multintlimits@{\intop\ifnum\intno@=\z@\intdots@
  \else\intkern@\fi
    \ifnum\intno@>\tw@\intop\intkern@\fi
     \ifnum\intno@>\thr@@\intop\intkern@\fi\intop}
\def\multint@{\int\ifnum\intno@=\z@\intdots@\else\intkern@\fi
   \ifnum\intno@>\tw@\int\intkern@\fi
    \ifnum\intno@>\thr@@\int\intkern@\fi\int}
\def\ints@@{\iflimtoken@\def\ints@@@{\iflimits@
   \negintic@\mathop{\intic@\multintlimits@}\limits\else
    \multint@\nolimits\fi\eat@}\else
     \def\ints@@@{\multint@\nolimits}\fi\ints@@@}
\def\Sb{_\bgroup\vspace@
        \baselineskip=\fontdimen10 \scriptfont\tw@
        \advance\baselineskip by \fontdimen12 \scriptfont\tw@
        \lineskip=\thr@@\fontdimen8 \scriptfont\thr@@
        \lineskiplimit=\thr@@\fontdimen8 \scriptfont\thr@@
        \Let@\vbox\bgroup\halign\bgroup \hfil$\scriptstyle
            {##}$\hfil\cr}
\def\endSb{\crcr\egroup\egroup\egroup}
\def\Sp{^\bgroup\vspace@
        \baselineskip=\fontdimen10 \scriptfont\tw@
        \advance\baselineskip by \fontdimen12 \scriptfont\tw@
        \lineskip=\thr@@\fontdimen8 \scriptfont\thr@@
        \lineskiplimit=\thr@@\fontdimen8 \scriptfont\thr@@
        \Let@\vbox\bgroup\halign\bgroup \hfil$\scriptstyle
            {##}$\hfil\cr}
\def\endSp{\crcr\egroup\egroup\egroup}
\def\Let@{\relax\iffalse{\fi\let\\=\cr\iffalse}\fi}
\def\vspace@{\def\vspace##1{\noalign{\vskip##1 }}}
\def\aligned{\,\vcenter\bgroup\vspace@\Let@\openup\jot\m@th\ialign
  \bgroup \strut\hfil$\displaystyle{##}$&$\displaystyle{{}##}$\hfil\crcr}
\def\endaligned{\crcr\egroup\egroup}
\def\matrix{\,\vcenter\bgroup\Let@\vspace@
    \normalbaselines
  \m@th\ialign\bgroup\hfil$##$\hfil&&\quad\hfil$##$\hfil\crcr
    \mathstrut\crcr\noalign{\kern-\baselineskip}}
\def\endmatrix{\crcr\mathstrut\crcr\noalign{\kern-\baselineskip}\egroup
                \egroup\,}
\newtoks\hashtoks@
\hashtoks@={#}
\def\format{\crcr\egroup\iffalse{\fi\ifnum`}=0 \fi\format@}
\def\format@#1\\{\def\preamble@{#1}%
  \def\c{\hfil$\the\hashtoks@$\hfil}%
  \def\r{\hfil$\the\hashtoks@$}%
  \def\l{$\the\hashtoks@$\hfil}%
  \setbox\z@=\hbox{\xdef\Preamble@{\preamble@}}\ifnum`{=0 \fi\iffalse}\fi
   \ialign\bgroup\span\Preamble@\crcr}

\def\cases{\left\{\,\vcenter\bgroup\vspace@
     \normalbaselines\openup\jot\m@th
       \Let@\ialign\bgroup$##$\hfil&\quad$##$\hfil\crcr
      \mathstrut\crcr\noalign{\kern-\baselineskip}}
\def\endcases{\endmatrix\right.}
\newif\iftagsleft@
\tagsleft@true
\def\TagsOnRight{\global\tagsleft@false}
\def\tag#1$${\iftagsleft@\leqno\else\eqno\fi
 \hbox{\def\pagebreak{\global\postdisplaypenalty-\@M}%
 \def\nopagebreak{\global\postdisplaypenalty\@M}\rm(#1\unskip)}%
  $$\postdisplaypenalty\z@\ignorespaces}
\interdisplaylinepenalty=\@M
\def\allowdisplaybreak@{\def\allowdisplaybreak{\noalign{\allowbreak}}}
\def\displaybreak@{\def\displaybreak{\noalign{\break}}}
\def\align#1\endalign{\def\tag{&}\vspace@\allowdisplaybreak@\displaybreak@
  \iftagsleft@\lalign@#1\endalign\else
   \ralign@#1\endalign\fi}
\def\ralign@#1\endalign{\displ@y\Let@\tabskip\centering\halign to\displaywidth
     {\hfil$\displaystyle{##}$\tabskip=\z@&$\displaystyle{{}##}$\hfil
       \tabskip=\centering&\llap{\hbox{(\rm##\unskip)}}\tabskip\z@\crcr
             #1\crcr}}
\def\lalign@
 #1\endalign{\displ@y\Let@\tabskip\centering\halign to \displaywidth
   {\hfil$\displaystyle{##}$\tabskip=\z@&$\displaystyle{{}##}$\hfil
   \tabskip=\centering&\kern-\displaywidth
        \rlap{\hbox{(\rm##\unskip)}}\tabskip=\displaywidth\crcr
               #1\crcr}}
\def\overrightarrow{\mathpalette\overrightarrow@}
\def\overrightarrow@#1#2{\vbox{\ialign{$##$\cr
    #1{-}\mkern-6mu\cleaders\hbox{$#1\mkern-2mu{-}\mkern-2mu$}\hfill
     \mkern-6mu{\to}\cr
     \noalign{\kern -1\p@\nointerlineskip}
     \hfil#1#2\hfil\cr}}}
\def\overleftarrow{\mathpalette\overleftarrow@}
\def\overleftarrow@#1#2{\vbox{\ialign{$##$\cr
     #1{\leftarrow}\mkern-6mu\cleaders\hbox{$#1\mkern-2mu{-}\mkern-2mu$}\hfill
      \mkern-6mu{-}\cr
     \noalign{\kern -1\p@\nointerlineskip}
     \hfil#1#2\hfil\cr}}}
\def\overleftrightarrow{\mathpalette\overleftrightarrow@}
\def\overleftrightarrow@#1#2{\vbox{\ialign{$##$\cr
     #1{\leftarrow}\mkern-6mu\cleaders\hbox{$#1\mkern-2mu{-}\mkern-2mu$}\hfill
       \mkern-6mu{\to}\cr
    \noalign{\kern -1\p@\nointerlineskip}
      \hfil#1#2\hfil\cr}}}
\def\underrightarrow{\mathpalette\underrightarrow@}
\def\underrightarrow@#1#2{\vtop{\ialign{$##$\cr
    \hfil#1#2\hfil\cr
     \noalign{\kern -1\p@\nointerlineskip}
    #1{-}\mkern-6mu\cleaders\hbox{$#1\mkern-2mu{-}\mkern-2mu$}\hfill
     \mkern-6mu{\to}\cr}}}
\def\underleftarrow{\mathpalette\underleftarrow@}
\def\underleftarrow@#1#2{\vtop{\ialign{$##$\cr
     \hfil#1#2\hfil\cr
     \noalign{\kern -1\p@\nointerlineskip}
     #1{\leftarrow}\mkern-6mu\cleaders\hbox{$#1\mkern-2mu{-}\mkern-2mu$}\hfill
      \mkern-6mu{-}\cr}}}
\def\underleftrightarrow{\mathpalette\underleftrightarrow@}
\def\underleftrightarrow@#1#2{\vtop{\ialign{$##$\cr
      \hfil#1#2\hfil\cr
    \noalign{\kern -1\p@\nointerlineskip}
     #1{\leftarrow}\mkern-6mu\cleaders\hbox{$#1\mkern-2mu{-}\mkern-2mu$}\hfill
       \mkern-6mu{\to}\cr}}}
\def\sqrt#1{\radical"270370 {#1}}
\def\dots{\relax\ifmmode\let\next=\ldots\else\let\next=\tdots@\fi\next}
\def\tdots@{\unskip\ \tdots@@}
\def\tdots@@{\futurelet\next\tdots@@@}
\def\tdots@@@{$\mathinner{\ldotp\ldotp\ldotp}\,
   \ifx\next,$\else
   \ifx\next.\,$\else
   \ifx\next;\,$\else
   \ifx\next:\,$\else
   \ifx\next?\,$\else
   \ifx\next!\,$\else
   $ \fi\fi\fi\fi\fi\fi}
\def\text{\relax\ifmmode\let\next=\text@\else\let\next=\text@@\fi\next}
\def\text@@#1{\hbox{#1}}
\def\text@#1{\mathchoice
 {\hbox{\everymath{\displaystyle}\def\textfonti{\the\textfont1 }%
    \def\textfontii{\the\textfont2 }\textdef@@ T#1}}
 {\hbox{\everymath{\textstyle}\def\textfonti{\the\textfont1 }%
    \def\textfontii{\the\textfont2 }\textdef@@ T#1}}
 {\hbox{\everymath{\scriptstyle}\def\textfonti{\the\scriptfont1 }%
   \def\textfontii{\the\scriptfont2 }\textdef@@ S\rm#1}}
 {\hbox{\everymath{\scriptscriptstyle}\def\textfonti{\the\scriptscriptfont1 }%
   \def\textfontii{\the\scriptscriptfont2 }\textdef@@ s\rm#1}}}
\def\textdef@@#1{\textdef@#1\rm \textdef@#1\bf
   \textdef@#1\sl \textdef@#1\it}

\def\textdef@#1#2{\def\next{\csname\expandafter\eat@\string#2fam\endcsname}%
\if S#1\edef#2{\the\scriptfont\next\relax}%
 \else\if s#1\edef#2{\the\scriptscriptfont\next\relax}%
 \else\edef#2{\the\textfont\next\relax}\fi\fi}
\scriptfont\itfam=\tenit \scriptscriptfont\itfam=\tenit
\scriptfont\slfam=\tensl \scriptscriptfont\slfam=\tensl
\mathcode`\0="0030
\mathcode`\1="0031
\mathcode`\2="0032
\mathcode`\3="0033
\mathcode`\4="0034
\mathcode`\5="0035
\mathcode`\6="0036
\mathcode`\7="0037
\mathcode`\8="0038
\mathcode`\9="0039
\def\Cal{\relax\ifmmode\let\next=\Cal@\else
     \def\next{\errmessage{Use \string\Cal\space only in math mode}}\fi\next}
\def\Cal@#1{{\fam2 #1}}
\def\bold{\relax\ifmmode\let\next=\bold@\else
   \def\next{\errmessage{Use \string\bold\space only in math
      mode}}\fi\next}\def\bold@#1{{\fam\bffam #1}}
\mathchardef\Gamma="0000
\mathchardef\Delta="0001
\mathchardef\Theta="0002
\mathchardef\Lambda="0003
\mathchardef\Xi="0004
\mathchardef\Pi="0005
\mathchardef\Sigma="0006
\mathchardef\Upsilon="0007
\mathchardef\Phi="0008
\mathchardef\Psi="0009
\mathchardef\Omega="000A
\mathchardef\varGamma="0100
\mathchardef\varDelta="0101
\mathchardef\varTheta="0102
\mathchardef\varLambda="0103
\mathchardef\varXi="0104
\mathchardef\varPi="0105
\mathchardef\varSigma="0106
\mathchardef\varUpsilon="0107
\mathchardef\varPhi="0108
\mathchardef\varPsi="0109
\mathchardef\varOmega="010A
\font\dummyft@=dummy
\fontdimen1 \dummyft@=\z@
\fontdimen2 \dummyft@=\z@
\fontdimen3 \dummyft@=\z@
\fontdimen4 \dummyft@=\z@
\fontdimen5 \dummyft@=\z@
\fontdimen6 \dummyft@=\z@
\fontdimen7 \dummyft@=\z@
\fontdimen8 \dummyft@=\z@
\fontdimen9 \dummyft@=\z@
\fontdimen10 \dummyft@=\z@
\fontdimen11 \dummyft@=\z@
\fontdimen12 \dummyft@=\z@
\fontdimen13 \dummyft@=\z@
\fontdimen14 \dummyft@=\z@
\fontdimen15 \dummyft@=\z@
\fontdimen16 \dummyft@=\z@
\fontdimen17 \dummyft@=\z@
\fontdimen18 \dummyft@=\z@
\fontdimen19 \dummyft@=\z@
\fontdimen20 \dummyft@=\z@
\fontdimen21 \dummyft@=\z@
\fontdimen22 \dummyft@=\z@
\def\fontlist@{\\{\tenrm}\\{\sevenrm}\\{\fiverm}\\{\teni}\\{\seveni}%
 \\{\fivei}\\{\tensy}\\{\sevensy}\\{\fivesy}\\{\tenex}\\{\tenbf}\\{\sevenbf}%
 \\{\fivebf}\\{\tensl}\\{\tenit}\\{\tensmc}}
\def\dodummy@{{\def\\##1{\global\let##1=\dummyft@}\fontlist@}}
\newif\ifsyntax@
\newcount\countxviii@
\def\newtoks@{\alloc@5\toks\toksdef\@cclvi}
\def\nopages@{\output={\setbox\z@=\box\@cclv \deadcycles=\z@}\newtoks@\output}
\def\syntax{\syntax@true\dodummy@\countxviii@=\count18
\loop \ifnum\countxviii@ > \z@ \textfont\countxviii@=\dummyft@
   \scriptfont\countxviii@=\dummyft@ \scriptscriptfont\countxviii@=\dummyft@
     \advance\countxviii@ by-\@ne\repeat
\dummyft@\tracinglostchars=\z@
  \nopages@\frenchspacing\hbadness=\@M}
\def\magstep#1{\ifcase#1 1000\or
 1200\or 1440\or 1728\or 2074\or 2488\or 
 \errmessage{\string\magstep\space only works up to 5}\fi\relax}
{\lccode`\2=`\p \lccode`\3=`\t 
 \lowercase{\gdef\tru@#123{#1truept}}}

\def\scaletype#1{\mag=#1\relax
 \hsize=\expandafter\tru@\the\hsize
 \vsize=\expandafter\tru@\the\vsize
 \dimen\footins=\expandafter\tru@\the\dimen\footins}

\def\scalefont#1#2\andcallit#3{\edef\font@{\the\font}#1\font#3=
  \fontname\font\space scaled #2\relax\font@}
\def\Mag@#1#2{\ifdim#1<1pt\multiply#1 #2\relax\divide#1 1000 \else
  \ifdim#1<10pt\divide#1 10 \multiply#1 #2\relax\divide#1 100\else
  \divide#1 100 \multiply#1 #2\relax\divide#1 10 \fi\fi}
\def\scalelinespacing#1{\Mag@\baselineskip{#1}\Mag@\lineskip{#1}%
  \Mag@\lineskiplimit{#1}}
\def\wlog#1{\immediate\write-1{#1}}
\catcode`\@=\active

\input amssym.def
\nopagenumbers
\TagsOnRight
\scaletype{\magstep1}
\scalelinespacing{\magstep1}
\pagewidth{15truecm}
\pageheight{22truecm}
\headline={\rightline{\folio}}
\voffset=2\baselineskip

\define\re{\text{Re\,}}

\define\Ai{\text{Ai\,}}

\title Shape fluctuations and random matrices\endtitle
\bigskip
\bigskip
\author Kurt Johansson \endauthor
\vskip 10truecm
\noindent
Version: Revised August 19, 1999
\newline Address: \newline Department of Mathematics\newline
Royal Institute of Technology\newline
S-100 44 Stockholm \newline
Sweden\newline
e-mail kurtj\@math.kth.se\newline
fax. +46 8 723 17 88

\newpage

\subheading{Abstract} We study a certain random growth model in two
dimensions closely related to the one-dimensional totally asymmetric
exclusion process. The results show that the shape fluctuations,
appropriately scaled, converges in distribution to the Tracy-Widom
largest eigenvalue distribution for the Gaussian Unitary Ensemble (GUE).
\bigskip\bigskip
\subheading{Mathematics Subject Classification}
\newline
Primary: 60C05, 60K35, 82C22, 82C41
\newline
Secondary: 05E10, 33C45, 15A52
\subheading{Keywords} Largest eigenvalue, random matrices, shape
fluctuations, asymmetric simple exclusion process, directed
first-passage percolation, Meixner polynomials

\newpage

\heading 1. Introduction and results\endheading
\medskip
The shape and height fluctuations in many 2-d random growth models
are expected to be of order
$N^\chi$, with $\chi=1/3$, if the mean of the linear size
of the shape or the height is of order $N$. See [KS] for a review and [NP] for
rigorous bounds on $\chi$ in first-passage percolation.

In this paper we will consider a specific model. It can be given
several probabilistic interpretations, as a randomly growing Young
diagram, a totally asymmetric one dimensional exclusion process, a
certain zero-temperature directed polymer in a random environment or
as a kind of first-passage site percolation model. The model has the
advantage that we can prove that $\chi=1/3$ and also compute the
asymptotic distribution of the appropriately rescaled random
variable. Interestingly, the limit distribution that occurs is the
same as that of the scaled largest eigenvalue of an $N\times N$ random
matrix from the Gaussian Unitary Ensemble (GUE) in the limit
$N\to\infty$. The model in this paper has many similarities with the
problem of the distribution of the length of the longest increasing
subsequence in a random permutation where the same limiting
distribution and $\chi=1/3$ was found in [BDJ]. 

To define the model 
let $w(i,j)$, $(i,j)\in\Bbb Z_+^2$, be independent
geometrically distributed random variables,
$$
\Bbb P[w(i,j)=k]=(1-q)q^k, \quad k\in\Bbb N,
$$
where $0<q<1$.
Let $\Pi_{M,N}$ be the
set of all up/right paths $\pi$ in $\Bbb Z_+^2$ from $(1,1)$ to
$(M,N)$, i.e. sequences $(i_k, j_k)$, $k=1,\dots ,M+N-1$, of sites in
$\Bbb Z_+^2$ such that $(i_1,j_1)=(1,1)$, $(i_{M+N-1},j_{M+N-1})=
(M,N)$ and $(i_{k+1},j_{k+1})-(i_k,j_k)=(1,0)$ or  $(0,1)$.
Define the random variable
$$
G(M,N)=\max_{\pi\in\Pi_{M,N}}\sum_{(i,j)\in\pi} w(i,j). \tag 1.1
$$
We also define the closely related random variable
$$
G^\ast(M,N)=\max_{\pi\in\Pi_{M,N}}\sum_{(i,j)\in\pi} w^\ast(i,j),
$$
where $w^\ast(i,j)=w(i,j)+1$, so that $\Bbb
P[w^\ast(i,j)=k]=(1-q)q^{k-1}$,
$k\ge 1$. Clearly,
$$
G^\ast(M,N)=G(M,N)+M+N-1,\tag 1.2
$$
since all paths have the same length.
Using this random variable we can define, for each $t\ge 0$, a random
subset of the first quadrant by
$$
A(t)=\{(M,N)\in\Bbb Z_+^2\,;\,G^\ast(M,N)\le t\}+[-1,0]^2.\tag 1.3
$$
From the definition of $G^\ast(M,N)$ and the fact that we consider up/right
paths it follows that $A(t)$ has the form
$$
\cup_{k=1}^r [k-1,k]\times[0,\lambda_k]
$$
for some integers $\lambda_1\ge\lambda_2\ge\dots\ge\lambda_r\ge1$, so
we can think of $A(t)$ as a Young diagram
$\lambda=(\lambda_1,\dots,\lambda_r)$. 
If we think of $t\in\Bbb N$ as a discrete time variable, $A(t)$ is a
randomly growing Young diagram. Let $\partial^\ast A(t)$ be those
unit cubes adjacent to $A(t)$ that can be added to $A(t)$ so that it
is still a Young diagram, i.e. each cube in $\partial^\ast A(t)$ must
have a cube in $A(t)$ or $\Bbb R^2\setminus [0,\infty)^2$
immediately below and to the left of it. The
fact that the $w^\ast(i,j)$\,:s are independent and geometrically
distributed random variables implies that $A(t+1)$ is obtained by
picking each cube in $\partial^\ast A(t)$ independently with
probability $p=1-q$ and adding those cubes that were picked to $A(t)$.
(Recall that $\Bbb P[w^\ast(i,j)=k+l|w^\ast(i,j)\ge k]=\Bbb
P[w(i,j)=l]$, $l\ge 0$, the lack of memory property.) The starting
configuration is $A(0)=\emptyset$ and $\partial^\ast A(0)=[0,1]^2$.
In this model $G^\ast(M,N)=k$ means that the box $[M-1,M]\times
[N-1,N]$ is added at time $k$. This growth model has been considered
in [JPS].

This randomly growing Young diagram can also, equivalently, be thought
of as a certain totally asymmetric exclusion process with discrete
time, compare [Ro] or [Li], p. 412. Let $C(t)=\partial
([0,\infty)^2\setminus A(t))$ and note that $C(t)$ consists of
vertical and horizontal line segments of length 1. To each vertical
line segment we associate a 1 and to each horizontal line segment a
0. If we read the numbers along $C(t)$, starting at infinity along the
y-axis and ending at infinity along the x-axis, we get an infinite
sequence $X(t)=(\dots,x_{-1}(t),x_{0}(t),x_{1}(0),x_{2}(0),\dots)$ of
0's and 1's, starting with infinitely many 1's and ending with
infinitely many 0's; we let $x_0$ be the last number we have before
passing through the line $x=y$. We can think of $X(t)$ as a
configuration of particles, where $x_k=1$ means that there is a
particle at $k$, whereas $x_k=0$ means that there is no particle at
$k$. The stochastic growth of $A(t)$ described above corresponds to
the following stochasic dynamics of the particle system. At time $t$
each particle independently moves to the right with probability $1-q$
provided there is no particle immediately to the right of
it. Otherwise it does not move. The starting configuration is $x_k(0)= 
1_{(-\infty,0]}(k)$. In this particle model $G^\ast(M,N)=k$ means that
the particle initially at position $-(N-1)$ has moved $M$ steps at
time $k$.

Our first result concerns the mean and large deviation properties of
$G(M,N)$. 
\proclaim{Theorem 1.1} For each $q\in(0,1)$ and $\gamma\ge 1$,
$$
\lim_{N\to\infty}\frac 1N\Bbb E[G([\gamma
N],N)]=\frac{(1+\sqrt{q\gamma})^2}{1-q}-1\doteq\omega(\gamma,q).\tag
1.4
$$
Also, $G([\gamma N],N)$ has the following large deviation
properties. There are functions $i(\epsilon )$ and $\ell(\epsilon)$
(which depend on $q$ and $\gamma$), so that, for any $\epsilon>0$,
$$
\lim_{N\to\infty}\frac 1{N^2}\log\Bbb P[G([\gamma N],N)\le 
N(\omega(\gamma,q)-\epsilon)]=-\ell(\epsilon)\tag 1.5
$$
and
$$
\lim_{N\to\infty}\frac 1{N}\log\Bbb P[G([\gamma N],N)\ge 
N(\omega(\gamma,q)+\epsilon)]=-i(\epsilon).\tag 1.6
$$
The functions $\ell (x)$ and $i(x)$ are $>0$ if $x>0$.
\endproclaim

Note that the existence of the limit (1.4) follows by a subadditivity
argument, so it is the explicit form of the constant that is
interesting. The large deviation result (1.6) has been obtained in
[Se2]. 
The theorem will be proved in section 2.

The theorem implies that $\frac 1t A(t)$ has an asymptotic shape $A_0$
as $t\to\infty$, in the sense that given any $\epsilon >0$ 
$$
(1-\epsilon)A_0\subseteq \frac 1tA(t)\subseteq (1+\epsilon)A_0
$$
for all sufficiently large $t$. 
It follows from the definition of $A(t)$, (1.3), and
theorem 1.1 that 
$$
A_0=\{(x,y)\in[0,\infty)^2\,;\,y+2\sqrt{qxy}+x\le 1-q\}.
$$
The boundary of $A_0$ consists of two line segments from the origin to
$(1-q,0)$ and $(0,1-q)$ and part of an ellipse that is tangent to the
$x$- and $y$-axes. 

We now want to understand the fluctuations of
$A(t)$ around its asymptotic shape $A_0$, i.e. the fluctuations of
$G([\gamma N],N)$ around $N\omega(\gamma,q)$. Before we can formulate
the result we need some preliminaries.
Let $\Ai (x)$ be the Airy function defined by
$$
\Ai(x)=\frac 1{2\pi}\int_{-\infty}^\infty e^{i(t+is)^3/3+ix(t+is)}dt,
$$
where $s>0$ is arbitrary.
Consider the {\it Airy kernel}
$$
A(x,y)=\frac{\Ai (x)\Ai'(y)-\Ai'(x)\Ai(y)}{x-y}, \tag 1.7
$$
as an integral kernel on $L^2[s,\infty)$. The Fredholm determinant
$$
F(s)=\det (I-A)\mid_{L^2[s,\infty)}=\sum_{k=0}^\infty\frac
{(-1)^k}{k!} \int_{[s,\infty)^k}\det(A(x_i,x_j))_{i,j=1}^k d^kx\tag
1.8
$$
is a distribution function. It is the distribution function of the
appropriately scaled largest eigenvalue of an  $N\times N$ random
matrix from the Gaussian Unitary Ensemble (GUE) in the limit
$N\to\infty$, the Tracy-Widom distribution, see [TW1]. 
The distribution function $F(s)$
can also be defined using a certain Painlev\'e II function,
$$
F(s)=\exp[-\int_s^\infty (x-s)u(x)^2dx],\tag 1.9
$$
where $u(x)$ is the unique solution of the Painlev\'e II equation
$$
u''=2u^3+xu,
$$
with the asymptotics $u(x)\sim \Ai (x)$ as $x\to\infty$. The fact that
the expressions (1.8) and (1.9) are equal is proved in [TW1].

\proclaim{Theorem 1.2} For each $q\in(0,1)$, $\gamma\ge 1$ and $s\in
\Bbb R$,
$$
\lim_{N\to\infty}\Bbb P[\frac{G([\gamma N],N)-N\omega(\gamma,q)}
{\sigma(\gamma,q)N^{1/3}}\le s]=F(s),\tag 1.10
$$
where
$$
\sigma(\gamma,q)=\frac{q^{1/6}\gamma^{-1/6}}{1-q}
(\sqrt{\gamma}+\sqrt{q})^{2/3}(1+\sqrt{q\gamma})^{2/3}.\tag
1.11
$$
\endproclaim
The theorem will be proved in section 3.
We have not proved convergence of the moments of the rescaled random
variable, see remark 2.5. This theorem should be compared with the
result obtained in [BDJ], that if $\ell_N(\sigma)$ is the length of a
longest increasing subsequence in a random permutation $\sigma\in S_N$
(all $N!$ permutations have the same probability), then
$$
\lim_{N\to\infty}\Bbb P[(\sqrt{N})^{-1/3}(\ell_N(\sigma)-2\sqrt{N})\le
s]=F(s).\tag 1.12
$$
Note that in both cases we have the same exponent $1/3$, the standard
deviation is $\sim (\text{mean})^{1/3}$

The proofs of theorems 1.1 and 1.2 are based on the following result
which will be proved in section 2.
\proclaim{Proposition 1.3} For any $M\ge N\ge 1$,
$$
\Bbb P[G(M,N)\le t]=\frac 1{\Cal Z_{M,N}}\sum\Sb h\in\Bbb
N^N\\\max\{h_i\} \le t+N-1\endSb \prod_{1\le i<j\le
N}(h_i-h_j)^2\prod_{i=1}^N \binom{h_i+M-N}{h_i} q^{h_i},\tag 1.13
$$
where $\Cal Z_{M,N}$ is the normalization constant (partition
function). 
\endproclaim

This remarkable formula should be compared with the formula for the
distribution function for the largest eigenvalue, $\lambda_{\max}$, of an
$N\times N$ random matrix from GUE,
$$
\Bbb P[\lambda_{\max}\le t]=\frac 1{Z_N}\int_{(-\infty, t]^N}\prod_{1\le
i<j\le N}(x_i-x_j)^2\prod_{j=1}^Ne^{-2Nx_j^2}d^Nx.\tag 1.14
$$
There is a clear similarity between the two expressions, so we can use
the ideas developed to investigate (1.14). Just as the right hand side
of (1.14) can be written as a Fredholm determinant, so can the right
hand side of (1.13). The kernel for (1.13) is the
\it Meixner kernel\rm,
$$
\align
&K_{M,N}(x,y)\\&=\frac{\kappa_{N-1}}{\kappa_N}
\frac{M_N(x)M_{N-1}(y)-M_{N-1}(x)M_N(y)}{x-y}
(w_K^q(x)w_K^q(y))^{1/2},\tag 1.15
\endalign
$$
where $M_N(x)=\kappa_N x^N+\dots$ are the normalized orthogonal
polynomials with respect to the discrete weight, $K=M-N+1$,
$$
w_K^q(x)=\binom{x+K-1}x q^x,\quad x\in \Bbb N.\tag 1.16
$$
This Meixner kernel also appears in the recent paper [BO].
The polynomial $M_N(x)$ is a multiple of the classical Meixner
polynomials $m_N^{K,q}(x)$. Using the explicit generating function for
the Meixner polynomials, see [Ch], the appropriate asymptotics of the
kernel (1.15) can be analyzed. This will be done in section 5.

Let $u(i,j)$, $(i,j)\in\Bbb Z_+^2$, be independent exponentially
distributed random variables with parameter 1. Let $H(M,N)$ be the
analogue of $G(M,N)$ for these random variables, i.e.
$$
H(M,N)=\max\{\sum_{(i,j)\in\pi}u(i,j)\,;\,\pi\in\Pi_{M,N}\}.\tag 1.17
$$
We can consider the related stochastically growing Young diagram and
totally asymmetric exclusion process just as in the geometric case,
where we now have continuous time. This simple exclusion process is
exactly the one considered by Rost, [Ro], see also [Li]. In this
process $X(t)=(\eta_k(t))_{k=-\infty}^\infty\in\{0,1\}^{\Bbb Z}$ the
initial configuration is $1_{(-\infty,0]}(k)$ and a particle
($\eta_k=1$) jumps with exponential rate to the right one step
provided there is no particle at $k+1$ ($\eta_{k+1}=0)$. By taking the
$q\to 1$ limit in (1.13) we obtain
\proclaim{Proposition 1.4} For any $M\ge N\ge 1$, $t\ge 0$,
$$
\Bbb P[H(M,N)\le t]=
\frac 1{Z_{M,N}'}\int_{[0, t]^N}\prod_{1\le
i<j\le N}(x_i-x_j)^2\prod_{j=1}^Nx_j^{M-N}e^{-x_j}d^Nx.\tag 1.18
$$
\endproclaim
\demo{Proof} If $X_L$ is geometrically distributed with parameter
$1-1/L$, then $L^{-1}X_L$ converges in distribution to an exponential
random variable with parameter 1. Since $G(M,N)$ is a continuous
function of the $w(i,j)$\,:s, proposition 1.3 gives
$$
\align
&\Bbb P[H(M,N)\le t]\\
&=\lim_{L\to\infty}\frac 1{\Cal Z_{M,N}}\sum_{(\ast)}
\prod_{1\le i<j\le
N}(h_i-h_j)^2\prod_{i=1}^N \binom{h_i+M-N}{h_i} (1-1/L)^{h_i}\\
&=\lim_{L\to\infty}\frac {L^{N^2}}{\Cal Z_{M,N}(M-N)!}\sum_{(\ast)}
\prod_{1\le i<j\le
N}(\frac{h_i-h_j}L)^2\prod_{i=1}^N 
e^{-\frac{h_i}L+o(\frac 1L)}\prod_{k=1}^{M-N}
(\frac{h_i+k}L)  \\
&=\frac 1{Z_{M,N}'}\int_{[0, t]^N}\prod_{1\le
i<j\le N}(x_i-x_j)^2\prod_{j=1}^Nx_j^{M-N}e^{-x_j}d^Nx,
\endalign
$$
where $(\ast)$ means summation over all $h\in\Bbb N^N$ such that 
$\max\{h_i\} \le [Lt]+N-1$.
\enddemo

\proclaim{Remark 1.5}\rm The right hand side in (1.18) is the probability that
the largest eigenvalue in the Laguerre ensemble is $\le t$.
It occurs in the following way. Let $A$ be an $N\times M$
rectangular matrix ($N\le M$) with entries that are complex Gaussian 
random variables with mean zero and variance 1/2. Then the right hand
side in (1.18) is the distribution function for the largest eigenvalue
of $AA^\ast$, see [Ja].\it
\endproclaim

\proclaim{Theorem 1.6} For each $\gamma\ge 1$,
$$
\lim_{N\to\infty}\frac 1N\Bbb E[H([\gamma
N],N)]=(1+\sqrt{\gamma})^2,\tag 1.19
$$
and there are functions $i_\ast(\epsilon)$ and $\ell_\ast(\epsilon)$
(which depend on $\gamma$), so that for any $\epsilon >0$,
$$
\lim_{N\to\infty}\frac 1{N^2}\log\Bbb P[H([\gamma
N],N)\le N((1+\sqrt{\gamma})^2-\epsilon)]=-\ell_\ast(\epsilon)\tag
1.20
$$
and
$$
\lim_{N\to\infty}\frac 1{N}\log\Bbb P[H([\gamma
N],N)\ge N((1+\sqrt{\gamma})^2+\epsilon)]=-i_\ast(\epsilon).\tag
1.21
$$
Furthermore, assume that $a_N=O(N^{1/3})$ as $N\to\infty$ and pick
$d_N$ so that $d_N-(1+1/\sqrt{\gamma})a_N=o(N^{1/3})$ as $N\to\infty$.
Then, for each $\gamma\ge 1$,
$$
\lim_{N\to\infty}\Bbb P[\frac{H(\gamma
N+a_N,N)-(1+\sqrt{\gamma})^2N-d_N}
{\gamma^{-1/6}(1+\sqrt{\gamma})^{4/3}N^{1/3}}\le s]=F(s).\tag 1.22
$$
\endproclaim
\demo{Proof} For the proof of (1.19) to (1.21) see remark 2.3. Write
$c=(1+\sqrt{\gamma})^2$ and
$\rho=\gamma^{-1/6}(1+\sqrt{\gamma})^{4/3}$. Then, by proposition 1.4,
$$\align
&\Bbb P[H(\gamma N+a_N,N)\le cN+d_N+\rho N^{1/3}s]\\&=
\frac 1{Z'_{\gamma N+a_N,N}}\int_{[0,cN+d_N+\rho N^{1/3}s]^N}\Delta
(x)^2\prod_{j=1}^Nx_j^{\alpha_N}e^{-x_j}d^Nx,
\endalign
$$
where $\Delta(x)=\prod_{1\le i<j\le N}(x_j-x_i)$ and
$\alpha_N=(\gamma-1)N+a_N$. By a standard argument, see [Me], ch. 5,
[TW3] or section
3, this equals the Fredholm determinant
$$
\sum_{k=0}^N\frac {(-1)^k}{k!}\int_{[s,\infty)^k}\det(\rho N^{1/3}
K_N^{\alpha_N}(cN+d_N+\rho N^{1/3}\xi_i,
cN+d_N+\rho N^{1/3}\xi_j))_{i,j=1}^kd^k\xi\tag 1.23
$$
where
$$
K_N^\alpha (x,y)=\frac{\kappa_{N-1}}{\kappa_N}\frac{\ell_N^\alpha(x)
\ell_{N-1}^\alpha(y)-\ell_N^\alpha(y)
\ell_{N-1}^\alpha(x)}{x-y}\left(x^\alpha e^{-x}y^\alpha
e^{-y}\right)^{1/2}.
$$
is the Laguerre kernel. Here, 
$$
\ell^\alpha_n(x)=\left(\frac{n!}{(\alpha+n)!}\right)^{1/2}(-1)^nL_n^\alpha
(x) =\kappa_n x^n+\dots
$$
are the normalized associated Laguerre polynomials,
$$
\int_0^\infty\ell^\alpha_n(x)\ell^\alpha_m(x)x^\alpha e^{-x}dx=\delta_{nm}.
$$
From asymptotic formulas for these polynomials it follows that
$$
\lim_{N\to\infty}K_N^{\alpha_N}(cN+d_N+\rho N^{1/3}\xi,
cN+d_N+\rho N^{1/3}\eta)=A(\xi,\eta).\tag 1.24
$$
This can be proved in the same way as the corresponding results for
Meixner polynomials, see sections 3 and 4, by using the integral
representation 
$$
L_n^\alpha(x)=\frac{e^x}{2\pi
i}\int_C\frac{e^{-xz}z^{n+\alpha}}{(z-1)^{n+1}}dz,
$$
where $C$ is a circle surrounding $z=1$. Using (1.23), (1.24) and some
estimates (compare lemma 3.1) we obtain
$$\align
&\lim_{N\to\infty}P[H(\gamma N+a_N,N)\le cN+d_N+\rho N^{1/3}s]\\
&=\sum_{k=0}^\infty\frac {(-1)^k}{k!}\int_{[s,\infty)^k}\det(A(\xi_i,\xi_j))
_{i,j=1}^kd^k\xi=F(s).
\endalign
$$
We will not present all the details since they are similar to the
proof of theorem 1.2.
\enddemo

Using this result we can get a fluctuation theorem for Rost's totally
asymmetric simple exclusion process defined above. The random variable
$H(N,M)$ is the first time at which the particle starting at $-(N-1)$ has
moved exactly $M$ steps to the right. If we define
$Y(k,t)=\sum_{j>k}\eta_j(t)$ 
to be the
number of particles to the right of $k$ at time $t$. Then $Y(k,t)>m$
means that the particle that starts at $-m$ has moved $\ge m+k+1$
steps at time $t$. Hence
$$
\Bbb P[Y(k,t)\le m]=1-\Bbb P[H(m+k+1,m+1)\le t].
$$
Using this relation and (1.19) to (1.21) we obtain the follwing result
first proved by Rost, [Ro],
$$
\frac 1tY([ut],t)\to\frac 14(1-u)^2 
$$
almost surely as $t\to\infty$, $|u|\le 1$. Now, using (1.22) it is
fairly straightforward to show the following result.
\proclaim{Corollary 1.7} For each $u\in [0,1)$,
$$
\lim_{t\to\infty}\Bbb P[Y([ut],t)\le\frac
t4 (1-u)^2+\frac{(1-u)^{2/3}}{(1+u)^{1/3}}\xi t^{1/3}]=1-F(-\xi).
$$
\endproclaim

\proclaim{Remark 1.8} \rm We can interpret theorems 1.1 and 1.2 (and
analogously theorem 1.6) as a result for a kind of zero-temperature
directed polymer or equivalently a directed first-passage site
percolation model in the following way.

Let $S_k$ be the simple random walk in $\Bbb Z$ starting at 0 at time
0 and ending at 0 at time $2N+2$. Denote the set of all possible paths
by $\Cal P_N$. Let $v(i,j)$, $(i,j)\in \Bbb Z^2$ be independent,
identically distributed random variables, and let $\beta>0$. On $\Cal
P_N$ we put the random path probability measure
$$
Q_N^\beta[S]=\frac 1{C_N^\beta}\exp(-\beta\sum_{k=1}^{2N}v(k,S_k)),
$$
$S\in\Cal P_N$, where $C_N^\beta$ is the normalization constant.
This measure describes a directed polymer ($S$) fixed at both
endpoints at inverse temperature $\beta$ in the random environment
given by the $v(i,j)$\,:s, see [Pi]. The {\it free energy} is $-\beta^{-1}\log
C_N^\beta$, and in the zero temperature limit $\beta\to\infty$ this
becomes
$$
F_N^{GS}=\min_{Z\in\Cal P_N}\sum_{k=1}^{2N}v(k,S_k),\tag 1.25
$$
the ground state energy. 
By rotating the coordinate system by the angle $-\pi/4$
it is seen that (1.25) can, equivalently, be thought of as a
first-passage time in a directed first passage site percolation model.
Let $u(i,j)$, $(i,j)\in\Bbb R_{+}^2$, be independent,
identically distributed random variables (with the same
distribution as the $v(i,j)$\,:s). Then
$F_N^{GS}$ has the same distribution as $F(N,N)$, where
$$
F(M,N)=\min_{\pi\in\Pi_{M,N}}\sum_{(i,j)\in\pi} u(i,j).
$$
(The $u(i,j)$\,:s are usually thought of as
passage times and $F(M,N)$ is the minimal flow time from $(1,1)$ to
$(M,N)$. Hence it is natural to assume that $u(i,j)\ge 0$, but this
will not be the case below.) 
We can define a random shape
$$
B(t)=\{(M,N)\in\Bbb Z_+^2;F(M,N)\le t\}+[-1,0]^2.
$$
Set $u(i,j)=\alpha-w(i,j)$, where $\alpha
>\alpha_{\text{min}} =(1-q)^{-1}(q+\sqrt{q})$ (this condition on
$\alpha$ ensures that $B(t)$ will grow); $w(i,j)$ are the
geometrically distributed random variables considered above.
Then clearly,
$$
F(M,N)=\alpha (M+N-1)-G(M,N),\tag 1.26
$$
Let
$\gamma\ge 1$, set $\hat x(\gamma)=(1+\gamma^2)^{-1/2}(\gamma, 1)$, a
unit
vector and $[n\hat x(\gamma)]=([N\gamma], N)$, ($[\cdot]$ the integer
part, where $N=[(1+\gamma^2)^{-1/2}n]$, so that $[n\hat x(\gamma)]$ is
a lattice site near $n\hat x(\gamma)$. Let $T_n(\gamma)$ be the first
time $s\ge 0$ for which $B(s)$ reaches $[n\hat x(\gamma)]$,
$$
T_n(\gamma)=\inf\{s\ge 0;[n\hat x(\gamma)]\in B(s)\}.
$$
Clearly, by the definition of $B(s)$ and equation (1.26),
$$
T_n(\gamma)=\alpha ([\gamma N]+N-1)-G([\gamma N],N),
$$
where $N=[(1+\gamma^2)^{-1/2}n]$.

Theorem 1.1 implies
that for each $q\in (0,1)$ and $\gamma\ge 1$,
$$
\lim_{n\to\infty}\frac 1n\Bbb E[T_n(\gamma)]=\frac 1{\sqrt{1+\gamma^2}}[\alpha 
(\gamma+1)-\frac{(1+\sqrt{q\gamma})^2}{1-q}+1]\doteq\mu(\gamma).
$$
Also, $T_n(\gamma)$ has large deviation
properties similar to those for $G([\gamma N],N)$. Using this result
we can compute the asymptotic shape of $B(t)$. It follows from theorem
1.2 that
$$
\Bbb
P[\frac{T_n(\gamma)-n\mu(\gamma)}{(1+\gamma^2)^{-1/6}\rho(q,\gamma)
n^{1/3}}\le s]\to 1-F(-s),
$$
as $n\to\infty$.
\it\endproclaim

\proclaim{Conjecture 1.9}\rm Is the result for $G([\gamma N],N)$
limited to geometric and exponential random variables? Normally, we
expect limit laws for appropriately scaled random variables to be
independent of the details. It is therefore natural to conjecture that
if the $w(i,j)$\,:s are i. i. d. random variables with some suitable
asumptions on their distribution, then there are constants $a$ and $b$
so that $(G([\gamma N],N)-aN)/bN^{1/3}$ converges to a random variable
with distribution $F(s)$. By remark 1.8 this leads to a related
conjecture for directed first-passage site percolation.

\bigskip
\heading 2. The Coulomb gas\endheading
\medskip\flushpar
\subheading{2.1 Combinatorics}
\medskip
The key combinatorial ingredient is the Knuth correspondence
introduced in [Kn]. It generalizes the Schensted correspondence [Sc]
which is used in [BDJ]. Write $[N]=\{1,\dots,N\}$. Let $\Cal M_{M,N}$
denote the set of all $M\times N$ matrices $A=(a_{ij})$ with
non-negative integer elements, and let $\Cal M_{M,N}^k$ be the subset
of those matrices that satisfy $\sum_{i=1}^M\sum_{j=1}^N a_{ij}=k$. A
two-rowed array
$$
\sigma=\left(\matrix i_1 &\dots &i_k\\ j_1 &\dots &
j_k\endmatrix\right)
$$
is called a generalized permutation if the columns $\binom{i_r}{j_r}$
are lexicographically ordered, i.e. either $i_r<i_{r+1}$ or
$i_r=i_{r+1}$, $j_r\le j_{r+1}$. There is a one-to-one correspondence
between the set $\Cal S_{M,N}^k$ of all generalized permutations of
length $k$, where the elements in the upper row come from $[M]$ and
the elements in the lower row from $[N]$, and $\Cal M_{M,N}^k$ defined
by $\sigma\to f(\sigma)=A=(a_{ij})$, where
$$
a_{ij}=\#\text{times}\,\,\,\binom ij \,\,\,\text{occurs
in}\,\,\sigma.
$$
We say that $\binom{i_{r_1}}{j_{r_1}},\dots,\binom{i_{r_m}}{j_{r_m}}$,
$r_1<r_2<\dots<r_m$ is an increasing subsequence in $\sigma$ if
$j_1\le j_2\le \dots\le j_{r_m}$. Let $\ell (\sigma)$ denote the
length of a longest increasing subsequence in $\sigma$.
\proclaim{Example}\rm The generalized permutatation
$$
\left(\matrix
1&1&1&2&2&2&3&3&4&4\\1&2&2&2&2&2&1&2&1&3\endmatrix\right)
$$
corresponds to
$$
\left(\matrix 1&2&0\\0&3&0\\1&1&0\\1&0&1\endmatrix\right).
$$
A longest increasing subsequence is 1 2 2 2 2 2 2 3 so
$\ell(\sigma)=8$.
\it\endproclaim
Recall from section 1 that $\Pi_{M,N}$ denotes the set of all up/right
paths $\pi$ from $(1,1)$ to $(M,N)$ through the sites $(i,j)$ with
$1\le i\le M$, $1\le j\le N$.
\proclaim{Lemma 2.1} For each $A\in\Cal M_{M,N}^k$,
$$
\max\{\sum_{(i,j)\in\pi}
a_{ij}\,;\,\pi\in\Pi_{M,N}\}=\ell(f^{-1}(A)).\tag 2.1
$$
\endproclaim
\demo{Proof} This is clear from the definitions. That we go to the
right corresponds to the fact that $i_{r_1}\le\dots\le i_{r_m}$  and that
we go up corresponds to $j_{r_1}\le\dots\le j_{r_m}$ (the upper row
gives row indices whereas the lower row gives column indices in the
matrix).
\enddemo
Now, Knuth has defined a one-to-one mapping from the set $\Cal
S_{M,N}^k$ to pairs $(P,Q)$ of semi-standard Young tableaux of the
same shape $\lambda$, which is a partition of $k$, $\lambda\vdash k$,
where $P$ has elements in $[N]$ and $Q$ has elements in $[M]$. (More
information on Young tableaux can be found in [Sa] and [Fu].) This
correspondence has the property that if $\sigma\to (P,Q)$ and $P, Q$
have shape $\lambda$, then $\ell(\sigma)=$ the length of the first
row, $\lambda_1$, in $\lambda$. Consider $G(M,N)$ defined by
(1.1). The $M\times N$ matrix $W=(w(i,j))$ is a random element
in $\Cal M_{M,N}$. Let
$$
S(M,N)=\sum_{i=1}^M\sum_{j=1}^N w(i,j)
$$ 
and
$$
p_{M,N}(t)=\Bbb P[G(M,N)\le t].
$$
Then,
$$
p_{M,N}(t)=\sum_{k=0}^\infty\Bbb P[G(M,N)\le t|S(M,N)=k]\Bbb
P[S(M,N)=k].\tag 2.2
$$
For a fixed $A\in\Cal M_{M,N}^k$ we have
$$
\Bbb P[\{A\}]=\prod_{i,j}(1-q)q^{a_{ij}}=(1-q)^{MN}q^k,
$$
since $\sum_{i,j}a_{ij}=k$. We have proved
\proclaim{Lemma 2.2} The conditional probability $\Bbb
P[\cdot|S(M,N)=k]$ is the uniform distribution on $\Cal M_{M,N}^k$.
\endproclaim

This lemma is the reason that we choose the $w(i,j)$\,:s to be
independent and geometrically distributed. Note that
$$
\Bbb P[S(M,N)=k]=\#\Cal M_{M,N}^k(1-q)^{MN}q^k.\tag 2.3
$$
Let $L(\lambda,M,N)$ denote the number of pairs $(P,Q)$ of
semi-standard Young tableaux of shape $\lambda$, such that $P$ has
elements in $[N]$ and $Q$ has elements in $[M]$. Combining lemma 2.1,
lemma 2.2 and the Knuth correspondence we see that
$$
\Bbb P[G(M,N)\le t|S(M,N)=k]=\frac 1{\#\Cal
M_{M,N}^k}\sum_{\lambda\vdash k,\lambda_1\le t}L(\lambda,M,N).\tag 2.4
$$
To compute $L(\lambda,M,N)$ we use
\proclaim{Lemma 2.3} The number of semi-standard tableaux
of shape $\lambda$ and elements in $[N]$ equals
$$
\prod_{1\le i<j\le N}\frac{\lambda_i-\lambda_j+j-i}{j-i}.
$$
\endproclaim
\demo{Proof} We have two formulas for the Schur polynomial 
in $N$ variables associated with the
partition $\lambda$, [Sa], [Fu],
$$
s_\lambda(x)=\sum_T x^T=
\frac{\det(x_j^{\lambda_i+N-i})_{1\le i,j\le N}}
{\det(x_j^{N-i})_{1\le i,j\le N}},
$$
where the sum is over all semi-standard $\lambda$-tableaux $T$ with
elements in $[N]$ and $x^T=x_1^{m_1}\dots x_N^{m_N}$ with $m_j$ equal
to the number of times $j$ occurs in $T$.
Hence, evaluating the Vandermonde determinants,
$$
s_\lambda(1,x,\dots,x^{N-1})=x^r\prod_{1\le i<j\le N}\frac
{x^{\lambda_i-\lambda_j+j-i} -1}{x^{j-i}-1},
$$
where $r=\sum_{i=1}^N(i-1)\lambda_i$. The number of semi-standard
tableaux with elements in $[N]$ equals
$$
s_\lambda(1,1,\dots,1)=\lim_{x\to 1} s_\lambda(1,x,\dots,x^{N-1})=
\prod_{1\le i<j\le N}\frac
{\lambda_i-\lambda_j+j-i}{j-i}.
$$
This completes the proof of the lemma.
\enddemo
It follows from lemma 2.3 that
$$
L(\lambda,M,N)=
\prod_{1\le i<j\le M}\frac
{\lambda_i-\lambda_j+j-i}{j-i}
\prod_{1\le i<j\le N}\frac
{\lambda_i-\lambda_j+j-i}{j-i}.\tag 2.5
$$

We assume from now on that $M\ge N$, the other case is analogous by
symmetry. Since the numbers in the columns in $P$ and $Q$ are strictly
increasing we must have $\lambda_i=0$ if $N<i\le M$. Hence
$$
L(\lambda,M,N)=
\prod_{1\le i<j\le M}\left(\frac
{\lambda_i-\lambda_j+j-i}{j-i}\right)^2\prod_{i=1}^N\prod_{j=N+1}^M\left( 
\frac{\lambda_i+j-i}{j-i}\right).
$$
Let $h_j=\lambda_j+N-j$, $j=1,\dots,N$, so that $h_1=\lambda_1+N-1$,
$h_N=\lambda_N\ge 0$ and $h_1>h_2>\dots>h_N$. Then
$$\align
&L(\lambda,M,N)=
\prod_{1\le i<j\le N}\frac{(h_i-h_j)^2}{(j-i)^2}
\prod_{i=1}^N\prod_{j=N+1}^M\frac{h_i+j-N}{j-i}\\
&=\prod_{j=0}^{N-1}\frac 1{j!(M-N+j)!}
\prod_{1\le i<j\le
N}(h_i-h_j)^2\prod_{i=1}^N\frac{(h_i+M-N)!}{h_i!}.\tag 2.6
\endalign$$
The condition $\sum_{j=1}^N\lambda_j=k$ translates into
$\sum_{j=1}^Nh_j= k+N(N-1)/2$ and $\lambda_1\le t$ to $h_1\le t+N-1$.
By (2.2), (2.3) and (2.4) we have
$$
p_{M,N}(t)=\sum_{k=0}^\infty (1-q)^{MN}q^k\sum_{\lambda\vdash
k,\lambda_1\le k}L(\lambda,M,N)
$$
and inserting (2.6) yields
$$\align
&p_{M,N}(t)= \frac{(1-q)^{MN}}{N!}q^{-N(N-1)/2}\prod_{j=0}^{N-1}\frac
1{j!(M-N+j)!}\\&\times\sum_{k=0}^\infty\sum\Sb h\in\Bbb N^N\\\sum
h_i=k+N(N-1)/2\\ \max\{h_i\}\le t+N-1\endSb
\prod_{1\le i<j\le
N}(h_i-h_j)^2\prod_{i=1}^N\frac{(h_i+M-N)!}{h_i!}
q^{\sum_{i=1}^Nh_i}.
\endalign
$$
where we have used the symmetry under permutation of the
$h_i$\,:s. Summing over $k$ gives all the possible values of
$\sum h_i$, so we obtain
$$
p_{M,N}(t)=\frac 1{\Cal Z_{M,N}}\sum\Sb h\in\Bbb
N^N\\\max\{h_i\}\le t+N-1\endSb 
\prod_{1\le i<j\le
N}(h_i-h_j)^2\prod_{i=1}^N w_{M-N+1}^q(h_i).\tag 2.7
$$
where $w_K^q(x)$ is given by (1.16) and
$$
\Cal Z_{M,N}=q^{N(N-1)/2}(1-q)^{-MN}\prod_{j=0}^{N-1}j!(M-N+j)!.\tag
2.8
$$
This proves proposition 1.3.
\medskip
\subheading{2.2 The large deviation estimate}
\medskip
In order to investigate the location of the rightmost charge in (2.7)
and prove large deviation formulas we rescale the discrete Coulomb gas
(2.7). Let $M=[\gamma N]$, $\gamma\ge 1$ fixed, and $K=K(N)=[\gamma
N]-M+1$. Set $\Bbb A_N=\frac 1N\Bbb N$, $\Bbb A_N(s)=\{x\in\Bbb A_N;\,
x\le s\}$ and
$$
V^{\gamma,q}_N(t)=-\frac 1N\log w_{K(N)}^q(Nt),\quad t\ge 0.
$$
Using Stirling's formula we see that
$$
\lim_{N\to\infty}V^{\gamma,q}_N(t)=t\log\frac
1q-(t+\gamma-1)\log(t+\gamma-1)+t\log t+(\gamma-1)\log(\gamma-1)\doteq
V^{\gamma,q}(t)\tag 2.9
$$
uniformly on compact subsets of $[0,\infty)$. 
(We will often omit the superscripts $\gamma$ and $q$.) 
Rescaling the variables
in (2.7) by setting $h_i=Nx_i$, $x_i\in\Bbb A_N$ we see that (2.7) can
be written
$$
p_N(t)\doteq p_{M(N),N}(t)=\frac{Z_N(\frac tN+1-\frac
1N)}{Z_N}, \tag 2.10
$$
where
$$
Z_N(s)=\sum_{x\in\Bbb
A_N(s)^N}\Delta_N(x)^2\exp\bigl(-N\sum_{j=1}^N V_N(x_j) \bigr)\tag
2.11
$$
and $Z_N=Z_N(\infty)$.
Here $\Delta_N(x)=\prod_{1\le i<j\le N}(x_j-x_i)$ is the Vandermonde
determinant. 

When investigating the large deviation properties of $p_N(t)$ we may
just as well consider more general confining potentials $V_N$. Assume
that $V_N:[0,\infty)\to\Bbb R$, $N\ge 1$, satisfy

\item{(i)} $V_N$ is continuous, $N\ge 1$.
\item{(ii)} There are constants $\xi>0,T\ge 0$ and $N_0\ge 1$ such
that
$$
V_N(t)\ge (1+\xi)\log(t^2+1)\tag 2.12
$$
for $t\ge T$ and $N\ge N_0$.
\item{(iii)} $V_N(t)\to V(t)$ uniformly on compact subsets of
$[0,\infty)$.

Set for $x\in\Bbb A_N^M$ and $\beta>0$ 
$$
Q_{M,N}(x)=|\Delta_M(x)|^\beta\prod_{j=1}^M\exp(-\frac{\beta N}2
V_N(x_j)).
$$
(This $M$ is not the same as the previous $M$.) Define the partition
functions 
$$ 
Z_{M,N}(t)=\sum_{x\in\Bbb A_N(t)^M} Q_{M,N}(x),
$$
$Z_{M,N}=Z_{N,M}(\infty)$
and the probability measure
$$
P_{M,N}[B]=\frac 1{Z_{M,N}}\sum_{x\in B}Q_{M,N}(x),
$$
$B\subseteq\Bbb N^M$.
We are interested in the distribution of the position of
the rightmost charge, $\max_{1\le k\le M} x_k$. Its distribution
function is given by
$$
F_{M,N}(t)=P_{M,N}[\max x_k\le t]=\frac{Z_{M,N}(t)}{Z_{M,N}}.\tag 2.13
$$
(If $M=N$ we write $F_N(t)$.)

In order to formulate the large deviation results for $F_N(t)$ we need
some results from weighted potential theory, [ST]. The results we need
differ from the usual ones since we are interested in the continuum
limit of a discrete Coulomb gas, so that the particle density of the
rescaled gas is always $\le 1$. Hence, the equilibrium measures will
be absolutely continuous with a density $\phi$ satisfying $0\le\phi\le
1$. Let $\Cal A_s$ denote the set of all $\phi\in L^1[0,s)$ such that
$0\le\phi\le 1$ and $\int_0^s\phi=1$, $1\le s\le\infty$. Given
$V:[0,\infty)\to\Bbb R$, continuous and such that there is a
$\delta>0$ and  a $T\ge 0$ such that 
$$
V(t)\ge (1+\delta)\log(t^2+1)\tag 2.14
$$
for $t\ge T$, we set
$$
k_V(x,y)=\log|x-y|^{-1}+\frac 12 V(x)+\frac 12 V(y)
$$
and
$$
I_V[\phi]=\int_0^s\int_0^sk_V(x,y)\phi(x)\phi (y)dxdy,
$$
for $\phi\in\Cal A_s$.

The proof of the next proposition is similar to the corresponding
result in weighted potential theory. See [DS] and also [LL] where a very
similar problem is treated.
\proclaim{Proposition 2.1} For each $s\in[1,\infty]$ there is a unique
$\phi_V^s\in\Cal A_s$ such that
$$
\inf_{\phi\in\Cal A_s}I_V[\phi]=I_V[\phi_V^s]=F_V^s.
$$
The extremal function $\phi_V^s$ has compact support. (If $s=\infty$
we will drop the superscript.)
\endproclaim

Let $b_V=\sup(\text{supp\,}\phi_V)$ be the right endpoint of the
support of $\phi_V$. Set $J(t)=0$ for $t\le b_V$ and
$$
J(t)=\inf_{\tau\ge t}\int_0^\infty k_V(\tau,x)\phi_V(x)dx-F_V\tag
2.15
$$
for $t\ge b_V$. Also, set
$$
L(t)=\frac 12 (F_V^t-F_V)
$$
for $t\ge 1$. The next theorem gives the large deviations for the
distribution function $F_N(t)$ defined by (2.13)
\proclaim{Theorem 2.2} Assume that $V_N(t)$ satisfies the assumptions
(i) - (iii) above. Then
$$
\lim_{N\to\infty}\frac 1{N^2}\log F_N(t)=-\beta L(t)\tag 2.16
$$
for any $t\ge 1$ and $L(t)>0$ if $t<b_V$. Assume furthermore that
$J(t)>0$ for $t>b_V$. Then 
$$
\lim_{N\to\infty}\frac 1N\log (1-F_N(t))=-\beta J(t)\tag 2.17
$$
for all $t$.
\endproclaim

We postpone the proof to section 4.
\proclaim{Remark 2.3} \rm
The same result is true for a continuous Coulomb
gas on $\Bbb R$ with density
$$
\frac 1{Z_N^\beta}|\Delta_N(x)|^\beta\exp(-\frac{\beta N}2\sum_{j=1}^N
V(x_j)), \tag 2.18
$$
on $\Bbb R^N$, which occur in random matrix theory. The choice
$\beta=2$ and $V(t)=2t^2$ corresponds to the Gaussian Unitary Ensemble
(GUE), compare (1.14). We assume that $V$ is continuous and satisfies
(2.14). In this case $\Cal A_s$ is replaced by $\Cal M_1(s)$, the set
of all probability measures on $(-\infty,s)$, and $\phi_V(x)dx$ is
replaced by the equilibrium measure $d\mu_V(t)$, see [Jo]. The proof
is essentially the same. The formula (2.16) for certain $V$ 
is a consequence of the
result in [BG], see also [HP]. Also, (2.17) has been proved in the case
$V(t)=t^2/2$ in [BDG]. If we take (2.18) on $[0,\infty)^N$ with
$\beta=2$ and $V(t)=-(M/N-1)\log t+t$ we get the measure in (1.18),
and in this way we can prove (1.19) to (1.21).
\it\endproclaim
We can now apply theorem 2.2 to the model we are interested in. It is
straightforward to verify that $V_N^{\gamma, q}$ satisfies the
conditions (i) - (iii) with limiting external potential
$V^{\gamma,q}(t)$. Write $b_{V^{\gamma,q}}=b(\gamma, q)$. The
computation of $\phi_{V^{\gamma,q}}$ will be outlined in
section 6. We have
$$
b(\gamma,q)=\frac{(1+\sqrt{q\gamma})^2}{1-q}.
$$
If $\gamma\ge 1/q$, then
$$
\phi_{V^{\gamma,q}}(t)=v(\frac 2c(t-a)-1),\quad a\le t\le b,
$$
where $a=\frac{(1-\sqrt{q\gamma})^2}{1-q}$, $c=b(\gamma,q)-a$ and
$$
v(x)=\frac 1{2\pi}[\arctan (\frac{Dx+1}{\sqrt{1-x^2}\sqrt{D^2-1}})-
\arctan (\frac{Bx+1}{\sqrt{1-x^2}\sqrt{B^2-1}})],\tag 2.19
$$
$B=(\gamma+q)/2\sqrt{q\gamma}$, $D=(1+q\gamma)/2\sqrt{q\gamma}$. If
$\gamma<1/q$, then, 
$$
\phi_{V^{\gamma,q}}(t)=\cases
1,&\text{if $0\le t\le a$}\\
v(\frac 2c(t-a)-1),&\text{if $a\le t\le b$},
\endcases
$$
where
$$
v(x)=\frac 1{2\pi}[\pi-\arctan (\frac{Dx+1}{\sqrt{1-x^2}\sqrt{D^2-1}})
-\arctan (\frac{Bx+1}{\sqrt{1-x^2}\sqrt{B^2-1}})]\tag 2.20
$$
with $a,c,B,D$ as before.

We will not discuss the explicit form of the lower tail rate
function. The upper tail rate function is given by
$$
J(t)=\frac c{8\sqrt{q\gamma}}\int_1^x(x-y)[\frac{\gamma-q}{y+B}+
\frac{1-q\gamma}{y+D}]\frac{dy}{\sqrt{y^2-1}},\tag 2.21
$$
with $c,B,D$ as above and $x=2(t-a)/c-1$. Using this formula we can
show that (see section 6) there are constants $c_1>0$ and $c_2>0$ so
that
$$
J(b+\delta)\ge
\cases
c_1\delta^{3/2} &\text{if $0\le\delta\le 1$}\\
c_2\delta &\text{if $\delta\ge1$}
\endcases\tag 2.22
$$
and
$$
J(b+\delta)=\frac{2(1-q)^{3/2}\gamma^{1/4}}{3q^{1/4}(\sqrt{q}+\sqrt{\gamma}) 
(1+\sqrt{q\gamma})}\delta^{3/2}+O(\delta^{5/2}).\tag 2.23
$$
In particular $J(t)>0$ if $t>b(\gamma,q)$. 

From (2.10), (2.13) and theorem 2.2 we obtain
$$
\lim_{N\to\infty}\frac 1{N^2}\log p_N(Nt)=-2L(t+1)\tag 2.24
$$
and
$$
\lim_{N\to\infty}\frac 1{N}\log(1- p_N(Nt))=-2J(t+1)\tag 2.25
$$
for each $t\ge 0$. These formulas imply theorem 1.1 with 
$\ell(\epsilon)=2L(b_V-\epsilon)$ and
$i(\epsilon)=2J(b_V+\epsilon)$. By theorem 2.2 and (2.22) we have
$i(\epsilon)>0$ and $\ell(\epsilon)>0$ if $\epsilon>0$.

By a superadditivity argument, the limit (2.25) actually gives a large
deviation estimate for all $N$, compare [Se1].

\proclaim{Corollary 2.4} For all $t\ge 0$ and $N\ge 1$,
$$
1-p_N(Nt)\le\exp(-2NJ(t+1)).\tag 2.26
$$
\endproclaim
\demo{Proof} For $1\le M_1\le M_2$ and $1\le N_1\le N_2$ we let
$G[(M_1,N_1), (M_2,N_2)]$ denote the maximum of $\sum_{(i,j)\in\pi}
w(i,j)$ over all up/right paths from $(M_1,N_1)$ to $(M_2,N_2)$.
Note that if $1\le M_1< M_2$ and $1\le N_1< N_2$, then
\item{(i)} $G[(M_1+1,N_1+1), (M_2,N_2)]$ and $G[(1,1),
(M_2-M_1,N_2-N_1)]$ are identically distributed.
\item{(ii)} $G[(1,1),(M_1,N_1)]$ and $G[(M_1+1,N_1+1), (M_2,N_2)]$ are
independent.
\flushpar Since $[2\gamma N]\ge2[\gamma N]$, we have
\item{(iii)} $G[([\gamma N]+1, N+1), ([2\gamma N],2N)]\ge
G[([\gamma N]+1, N+1), (2[\gamma N],2N)]$.

Write $a_N(t)=1-p_N(Nt)=\Bbb P[G((1,1),([\gamma N],N))>Nt]$. 
Then, by (i) and (iii),
$$
a_N(t)\le\Bbb P[G(([\gamma N]+1, N+1), ([2\gamma
N],2N))>Nt]
$$
and hence, by (ii), $a_N(t)^2\le a_{2N}(t)$.
Repeated use of this inequality yields $N^{-1}\log
a_N(t)\le(2^kN)^{-1}\log a_{2^kN}(t)$ and by letting $k\to\infty$ and
using (2.25) we find $N^{-1}\log a_N(t)\le -2J(t+1)$.

\proclaim{Remark 2.5}\rm We cannot prove convergence of the moments of
the rescaled random variable in theorem 1.2 since we have no finite $N$
estimate
of $\Bbb P[G([\gamma N],N)-\omega N\le -sN^{1/3}]$ for $s>0$
large. This would require an estimate of the finite $N$ Fredholm
determinant. In
the other direction we can use the estimate in corollary 2.4. The same
remark applies to theorem 1.6.\it
\endproclaim
\proclaim{Remark 2.6} \rm In [BR] it is proved by Baik and Rains that
if we consider permutations with certain restrictions we can get the
Tracy-Widom distributions for GOE and GSE as limiting laws for longest
increasing and decreasing subsequences. By considering a restricted
geometry we can obtain the Tracy-Widom distribution for GOE, [TW2],
also in the present setting. Let $w(i,j)$ , $1\le i\le j$ be independent
geometrically distributed random variables, $\Bbb
P[w(i,j)=k]=(1-q)q^k$ for $1\le i<j$ and $\Bbb
P[w(i,i)=k]=(1-\sqrt{q})q^{k/2}$ for $i\ge 1$. Set $w(i,j)=w(j,i)$, if
$i>j\ge 1$, so that $A=(w(i,j))$ is a symmetric matrix. The Knuth
correspondence maps $A$ to a pair of semistandard Young tableaux
$(P,Q)$ with $Q=P$, i.e. $A$ maps to a single semistandard Young
tableaux, see [Kn] or [Fu]. Let $\Pi_{N,N}^{\text{sym}}$ be the set of
all up/right paths from $(1,1)$ to $(N,N)$ in $\{(i,j)\in\Bbb
Z_+^2\,;\, 1\le i\le j\}$, i. e. in a triangle, and set
$$
F(N)=\max\{\sum_{(i,j)\in\pi}w(i,j)\,;\,\pi\in
\Pi_{N,N}^{\text{sym}}\}.
$$
Now, we also have
$$
F(N)=\max\{\sum_{(i,j)\in\pi}w(i,j)\,;\,\pi\in\Pi_{N,N}\},
$$
which equals the length of the first row in $P$, because those parts
of a maximal path in $\Pi_{N,N}$ which goes below the diagonal can be
reflected in the diagonal to give a path in $\Pi_{N,N}^{\text{sym}}$
without changing the sum $\sum w(i,j)$ since $w(i,j)$ is symmetric.

The same argument as above now gives
$$
\Bbb P[F(N)\le t]=\frac 1{Z_N^{(1)}}\sum\Sb h\in\Bbb
N^N\\\max\{h_j\}\le t+N-1\endSb \prod_{1\le i<j\le
N}|h_i-h_j|\prod_{i=1}^Nq^{h_i/2}.
$$
This corresponds to $\beta = 1$, $\gamma =1$ in theorem 2.2. It should
be possible to
analyze the asymptotics in this case 
analogously to GOE, see [TW2], to show that we can
find constants $a$ and $b$ so that $\Bbb P[F(N)\le aN+sbN^{1/3}]$
conerges to $F_1(t)$, the Tracy-Widom distribution for GOE. 
However it is not immediate to generalize the techniques of [TW2],
so this remains to be done. Note that again we can take the limit $q\to
1$ to get the case of exponentially distributed random variables.

\enddemo
\bigskip
\heading 3. The Fredholm determinant\endheading
\medskip
From the identity (2.7) we have
$$
p_N(t)=\psi_N(t+N-1),\tag 3.1
$$
where
$$
\psi_N(s)=\Cal E_N[\prod_{j=1}^N(1-\chi_s(h_j))].\tag 3.2
$$
Here
$$
\Cal E_N[\cdot]=\frac 1{\Cal Z_{M(N),N}}\sum_{h\in\Bbb N^N}
(\cdot)\Delta_N(h)^2\prod_{j=1}^N w_{K(N)}^q(h_j),
$$
$K(N)=M(N)-N+1$, $M(N)=[\gamma N]$ and $\chi_s(t)$ is the indicator
function for the interval $(s,\infty)$. We will take $s$ in (3.2) to
be an integer. 

Let $M_j^{K,q}(x)$, $j=0,1,\dots$ be the normalized
orthogonal polynomials with respect to the weight $w_K^q(x)$ on $\Bbb
N$,
$$
\sum_{x=0}^\infty M_i^{k,q}(x)M_j^{K,q}(x)w_K^q(x)=\delta_{ij},
$$
and $M_j^{K,q}(x)=\kappa_jx^j+\dots$ with $\kappa_j>0$. Set
$$
K_N(x,y)=\sum_{j=0}^{N-1}M_j^{K,q}(x)M_j^{K,q}(y)w_K^q(x)^{1/2}
w_K^q(y)^{1/2},
$$
so that $K_N(x,y)$ is a reproducing kernel on $\ell^2(\Bbb
N)$. 

The polynomials $M_n^{K,q}$ are multiples of the standard Meixner
polynomials, [NSU], [Ch],
$$
M_n^{K,q}(x)=\frac{(-1)^n}{d_n}m_N^{K,q}(x),
$$
where
$$
d_n^2=\frac{n!(n+K-1)!}{(1-q)^Kq^n(K-1)!}.
$$
The leading coefficient in $m_n^{K,q}$ is $(\frac{q-1}q)^n$ and
consequently
$$
\kappa_n=\frac 1{d_n}\biggl(\frac{1-q}q\biggr)^n.
$$
The Meixner polynomials have the generating function, [Ch],
$$
\sum_{n=0}^\infty m_n^{K,q}(x)\frac{t^n}{n!}=(1-\frac
tq)^x(1-t)^{-x-K}.\tag 3.3
$$
The Christoffel-Darboux formula, [Sz], gives
$$
\align
&K_N(x,y)=\frac{\kappa_{N-1}}{\kappa_N}\frac{M_N(x)
M_{N-1}(y)-M_N(y)M_{N-1}(x)}{x-y}w_K^q(x)^{1/2}w_K^q(y)^{1/2}\\
&=-\frac{q}{(1-q)d_{N-1}^2}\frac {m_N(x)
m_{N-1}(y)-m_N(y)m_{N-1}(x)}{x-y}w_K^q(x)^{1/2}w_K^q(y)^{1/2},\tag 3.4
\endalign
$$
where we have omitted the upper indices.
Standard computations from random matrix theory, [Me], Ch. 5 and [TW2],
show that $\psi_N$ can be written as a Fredholm determinant,
$$
\psi_N(s)=\sum_{k=0}^N\frac {(-1)^k}{k!}\sum_{h\in\{s+1,s+2,\dots\}^k}
\det(K_N(h_i,h_j))_{1\le i,j\le k}.\tag 3.5
$$
The proof of theorem 1.2 is based on taking the appropriate limit in
(3.5).

The next lemma will allow us to compute the asymptotics of the right
hand side of (3.5).

\proclaim{Lemma 3.1} Let $b\ge 0$ be a constant and assume that
$\rho_N\to\infty$ as $N\to\infty$. Suppose furthermore that $K_N:
\Bbb N\times\Bbb B\to\Bbb R$, $N\ge 1$, satisfies the following
properties.
\item{(i)} Let $M_1>0$ be a given constant. There is a constant $C$
such that 
$$
\sum_{m=1}^\infty K_N(bN+\rho_N\tau+m,bN+\rho_N\tau+m)\le C\tag 3.6
$$
for all $N\ge 1$, $\tau\ge -M_1$.
\item{(ii)} Given $\epsilon >0$, there is an $L>0$ so that
$$
\sum_{m=1}^\infty K_N(bN+\rho_NL+m,bN+\rho_NL+m)\le \epsilon,\tag 3.7
$$
for all $N\ge 1$.
\item{(iii)} Let $M_0>0$ be a given constant. If $A(\xi,\eta)$ is the
Airy kernel defined by (1.7), then
$$
\lim_{N\to\infty}\rho_N
K_N(bN+\rho_N\xi,bN+\rho_N\eta)=A(\xi,\eta)\tag 3.8
$$
uniformly for $\xi,\eta\in [-M_0,M_0]$.
\item{(iv)} The matrix $(K_N(x_i,x_j))_{i,j=1}^k$ is positive definite
for any $x_i,x_j\in [0,\infty)$, $k\ge 1$

Then, for each fixed $t\in\Bbb R$,

$$
\lim_{N\to\infty}\sum_{k=0}^N\frac{(-1)^k}{k!}\sum_{h\in\Bbb N^k}
\det(K_N(bN+\rho_Nt+h_i,bN+\rho_Nt+h_j))_{i,j=1}^k
=F(t),\tag 3.9
$$
where $F(t)$ is given by (1.8).
\endproclaim

\demo{Proof} It follows from (iv) that
$$
|\det(K_N(x_i,x_j))_{1\le i,j\le k}|\le\prod_{j=1}^kK_N(x_j,x_j),\tag
3.10
$$
see for example [HJ]. Consequently,
$$
|\sum_{h\in\Bbb N^k}
\det(K_N(a_N+h_i,a_N+h_j))_{1\le i,j\le
k}|\le\biggl(\sum_{m=1}^\infty K_N(m,m)\biggr)^k.\tag 3.11
$$
where we have written $a_N=bN+\rho_N t$.

Choose $M_1$ so that $|t|\le M_1$.
Let $\epsilon >0$ be given. It follows from the estimates (3.6) and
(3.11) that we can choose $\ell$ so that
$$
|\sum_{k=\ell+1}^N\frac{(-1)^k}{k!}\sum_{h\in\Bbb N^k}
\det(K_N(a_N+h_i,a_N+h_j))_{i,j=1}
^{k}|\le\sum_{k=\ell+1}^\infty\frac{C^k}{k!}\le\epsilon,\tag 3.12
$$
for all $N\ge 1$. Choose $L_0$ so that (3.11) holds with $L=L_0-M_0$. 
Then, by the
estimates (3.6), (3.7) and (3.10), 
$$\align
&\left|\left(\sum_{h\in\Bbb N^k}-\sum_{h\in([L_0\rho_N]^c)^k}\right)
\det(K_N(a_N+h_i,a_N+h_j))_{1\le i,j\le
k}\right|\\
&\le\sum\Sb h\in\Bbb
N^k\\\text{some}\,\,h_j>L_0\rho_N\endSb\prod_{i=1}^k 
K_N(a_N+h_i,a_N+h_i)\\&\le\sum_{j=1}^k
\sum\Sb h\in\Bbb
N^k\\h_j>L_0\rho_N\endSb\prod_{i=1}^k 
K_N(a_N+h_i,a_N+h_i)\\
&\le k\biggl(\sum_{m=1}^\infty K_N(a_N+m,a_N+m)\biggr)^{k-1}
\biggl(\sum_{m=1}^\infty K_N(bN+L\rho_N+m,bN+L\rho_N+m)\biggr)\\&\le
kC^{k-1}\epsilon .\tag 3.13
\endalign
$$
Denote the Fredholm determinant in the right hand side of (3.9) by
$D_N(t)$. 
Inserting the estimates (3.12) and (3.13) into the formula (3.9) we
obtain
$$
\align
&\left|D_N(t)-\sum_{k=0}^\ell\frac{(-1)^k}{k!}\sum_{h\in
[L_0\rho_N]^k} 
\det(\Cal K_N(\sigma+\frac{h_i}{\rho_N},\sigma+\frac{h_j}
{\rho_N}))_{1\le i,j\le k}\frac 1{\rho_N^k}\right|\\
&\le\biggl(\sum_{k=0}^\ell\frac{kC^{k-1}}{k!}+1\biggr)\epsilon\le
(1+e^{C})\epsilon,\tag 3.14
\endalign
$$
where
$$
\Cal K_N(\xi,\eta)=\rho_N K_N(bN+\rho_N\xi,bN+\rho_N\eta).
$$
By assumption (iii), with $M_0=L_0+M_1$, 
we can chooose $N_0$ so that if $N\ge N_0$, then
$$
|\det(\Cal K_N(\sigma+\frac{x}{\rho_N},\sigma+\frac{y}
{\rho_N}))-\det(A(\sigma+\frac{x}{\rho_N},\sigma+\frac{y}
{\rho_N}))|\le\frac{\epsilon}{L_0^k}
$$
for all $x,y\in [L_0\rho_N]$. Thus,
$$\align
&\left|\sum_{k=0}^\ell\frac{(-1)^k}{k!}\sum_{h\in[L_0\rho_N]^k}\left[
\det(\Cal K_N(t+\frac{h_i}{\rho_N},t+\frac{h_j}
{\rho_N}))-\det(A(t+\frac{h_i}{\rho_N},t+\frac{h_j}
{\rho_N}))\right]\frac 1{\rho_N^k}\right|\\
&\le\sum_{k=0}^\ell\frac
1{k!}\left(\frac{L_0\rho_N+1}{L_0\rho_N}\right)^k\epsilon \le
C'\epsilon.\tag 3.15
\endalign
$$
Combining the estimates (3.14) and (3.15) we find
$$
\left|D_N(t)-
\sum_{k=0}^\ell\frac{(-1)^k}{k!}\sum_{h\in[L_0\rho_N]^k}
\det(A(\sigma+\frac{h_i}{\rho_N},\sigma+\frac{h_j}
{\rho_N}))_{i,j=1}^k\frac 1{\rho_N^k}\right|\le C''\epsilon.\tag 3.16
$$
The Airy kernel can be written, [TW1],
$$
A(x,y)=\int_0^\infty\Ai (x+s)\Ai(y+s)ds.\tag 3.17
$$
Using the formula, see for example [H\"o], p. 214,
$$
\Ai(x)=e^{-\frac 23x^{3/2}}\frac 1{2\pi}\int_{-\infty}^\infty
e^{-\xi^2\sqrt{x}+i\xi^3/3}d\xi,
$$
valid for $x>0$, we see that
$$
|\Ai(x)|\le\frac 1{2\sqrt{\pi}x^{1/4}}e^{-\frac 23x^{3/2}}, \quad x>0.
$$
This estimate can be used to show that the Airy kernel satisfies (i)
and (ii) above. Since the matrix $(A(\xi_i,\xi_j))_{1\le i,j\le
k}$ is positive definite, we can use the same argument as above to
show that
$$
\left|\left(\sum_{k=0}^\infty
\int_{[t,\infty)^k}
-
\sum_{k=0}^\ell\int_{[t,L_0]^k}\right)\frac{(-1)^k}{k!}
\det(A(\xi_i,\xi_j))_{i,j=1}^kd^k\xi\right|\le\epsilon\tag 3.18
$$
provided $\ell$ and $L_0$ are sufficiently large. From (3.17) we see
that choosing $N_1\ge N_0$ large enough we have
$$
\left|D_N(t)-\sum_{k=0}^\ell\frac{(-1)^k}{k!}\int_{[t,L_0]^k}
\det(A(\xi_i,\xi_j))_{1\le i,j\le k}d^k\xi\right|\le C'''\epsilon\tag
3.19 
$$
for all $N\ge N_1$. If we combine the estimates (3.18) and (3.19) we
have proved the lemma.
\enddemo

To apply this lemma to the Meixner kernel (3.4) we need

\proclaim{Lemma 3.2} The Meixner kernel satisfies the properties (i)
to (iv) in lemma 3.1 with $b=b(\gamma,q)$ as before and
$\rho_N=\sigma N^{1/3}$, where $\sigma$ is given by (1.11).
\endproclaim

This lemma will be proved in section 5. We can now combine (3.1),
(3.5) and (3.9) to get
$$
\lim_{N\to\infty} p_N((b-1)N+\sigma N^{1/3}t)=F(t),\tag 3.20
$$
which is (1.10) and theorem 1.2 is proved. 

\bigskip
\heading 4. Proof of the large deviation theorem\endheading
\medskip
In this section we will prove theorem 2.2.
Set
$$
K_{N,V}=\sum_{1\le i\neq j\le N}k_V(x_i,x_j).
$$
By adding a constant $C$ to $V_N$, which does not alter the problem
we can, by assumption (ii) on $V_N$, assume that
$$
V_N(t)-\log(t^2+1)\ge\xi\log(t^2+1)\tag 4.1
$$
for all $t\ge 0$. Since $|t-s|^2\le (t^2+1)(s^2+1)$, this implies
$$
-K_{M,V_N}(x)\le -\xi (M-1)\sum_{j=1}^{M-1}\log(1+x_j^2)\tag 4.2
$$
for all $x\in [0,\infty)^M$.
Note that
$$
\sum_{1\le j\neq k\le N-1}\log|x_j-x_k|-N\sum_{j=1}^{N-1}V_N(x_j)=
-K_{N-1,V_N}(x)-\sum_{j=1}^{N-1}V_N(x_j).\tag 4.3
$$
The next lemma is analogous to lemma 4.2 in [Jo].
\proclaim{Lemma 4.1} Let $\{s_N\}$ be a sequence in $[0,\infty)$ such
that $s_N\to s>0$ as $N\to\infty$, or $s_N\equiv\infty$. Set, for
a given $\alpha >0$,
$$
\Omega_{N,\alpha}(s)=\{x\in\Bbb A_N(s)^{N-1}\,;\,\frac
1{N^2}K_{N-1,V_N}(x)\le F_V^\sigma+\alpha\}.
$$
Let $0\le\lambda\le 1$ and let $\sigma_N\in \Bbb A_N$, $N\ge 1$, be a
sequence converging to $\sigma>0$. Define a probability measure on $\Bbb
A_N(s_N)^{N-1}$ by
$$
P_{N-1,N}^{\lambda,\sigma_N}(\Omega;s_N)= \frac
1{Z_{N-1,N}^{\lambda,\sigma_N}
(s_N)}\sum_{x\in\Omega}\prod_{j=1}^{N-1} |\sigma_N-x_j|^{\lambda\beta}
Q_{N-1,N}(x), \tag 4.4
$$
where $Z_{N-1,N}^{\lambda,\sigma_N}(s_N)$ is a normalization constant.
($E_{N-1,N}^{\lambda,\sigma_N}[\cdot;s_N]$ denotes the corresponding
expectation and if $s_N\equiv\infty$ or $\lambda=\sigma_N=0$ we omit
them in the notation.) Fix $\eta>0$. Then there is an $N_1$ such that
for all $a\ge 0$ and $N\ge N_1$,
$$
P_{N-1,N}^{\lambda,\sigma_N}(\Omega_{N,\eta+a}(s_N)^c;s_N)\le e^{-\frac{\beta}4
aN^2}. \tag 4.5
$$
\endproclaim
\demo{Proof} We first prove the following claim.
\proclaim{Claim 4.2} Let $\sigma_N\in\Bbb A_N$, $\sigma_N\to\sigma$ as
$N\to\infty$ and $s\in
(0,\infty]$. For each $N\ge 2$ we can choose $(x_1^N,\dots
,x_{N-1}^N)\in \Bbb A_N(s)^{N-1}$ so that
$$
\frac 1{N^2}\sum_{1\le j\neq k\le N-1}\log|x_j^N-x_k^N|^{-1}+\frac 1N
\sum_{j=1}^{N-1} V_N(x_j^N)-\frac 1{N^2}\sum_{j=1}^{N-1}\log
|\sigma_N-x_j^N| \to F_V^s\tag 4.6
$$
as $N\to\infty$.
\endproclaim

To see this set
$$ 
y_k^N=\max\{\frac jN\,;\,\text{$j\in\Bbb N$ and $\int_0^{j/N}\phi_V^s
(t)dt<\frac kN$}\}.
$$
If $y_k^N\neq\sigma_N$ for $k=1,\dots,N-1$, we set $x_k^N=y_k^N$. If
$y_{k_0}^N=\sigma_N$, we set $x_k^N=y_k^N$ for $k<k_0$ and
$x_k^N=y_k^N+1/N$ for $k=k_0,\dots, N-1$. Using the fact that
$0\le\phi_V^s\le 1$ it is not difficult to see that
$x_1^N<x_2^N<\dots<x_{N-1}^N\le L$ for all $N$ and some fixed
$L$. Furthermore
$$
\frac 1{N-1}\sum_{k=1}^{N-1}\delta_{x_k^N}\to\phi_V^s(x)dx\tag 4.7
$$
weakly as $N\to\infty$. The property (iii) in the assumptions on $V_N$
implies
$$
\frac 1N\sum_{j=1}^{N-1}V_N(x_j^N)\to\int_0^\infty
V(t)\phi_V^s(t)dt.\tag 4.8
$$
Clearly,
$$
\frac 1{N^2}\sum_{j=1}^{N-1}\log|\sigma_N-x_j^N|^{-1}\le\frac 2{N^2}
\sum_{j=1}^{N-1}\log\frac Nj=\frac
2{N^2}\log\frac{N^{N-1}}{(N-1)!},\tag 4.9
$$
which $\to 0$ as $N\to\infty$. 
Also, since $\sigma_N\to\sigma$ and the $x_j^N$ belong to a bounded
set, we get a bound in the other direction which goes to 0 as $N\to\infty$.
Given $M\ge 1$, set
$f_M(t)=\min\{\log|t|^{-1},\log M\}$. Write 
$$
\align
&\frac 1{N^2}\sum_{j\neq k}\log|x_j^N-x_k^N|^{-1}=
\frac 1{N^2}\sum_{j\neq k}f_M(x_j^N-x_k^N)\\
&+\frac 1{N^2}\sum \Sb j\neq k\\|x_j^N-x_k^N|<1/M\endSb
(\log|x_j^N-x_k^N|^{-1}-f_M(x_j^N-x_k^N)).\tag 4.10
\endalign
$$
The absolute value of the second sum in the right hand side of (4.10)
is
$$
\le\frac 1{N^2}\sum\Sb 1\le|j-k|\le N/M\\|j|,|k|\le LN\endSb\log|\frac
N{j-k}| \le C\frac{\log M}{M}.
$$
Thus, using the weak convergence (4.7) and then letting $M\to\infty$
we obtain
$$
\lim_{N\to\infty}\frac 1{N^2}\sum_{j\neq k}\log|x_j^N-x_k^N|^{-1}
=\int_0^\infty\int_0^\infty\log|x-y|^{-1}\phi_V^s(x)\phi_V^s(y)dxdy,
$$
which together with (4.8) and (4.9) proves the claim.

We turn now to the proof of lemma 4.1. Let $\epsilon>0$ be given. We
want to estimate $Z_{N-1,N}^{\lambda,\sigma_N}$ from below. Choose
$N_0$ so that $s_N\ge s-\epsilon$ if $N\ge N_0$. Then
$$
Z_{N-1,N}^{\lambda,\sigma_N}(s_N)\ge
Z_{N-1,N}^{\lambda,\sigma_N}(s-\epsilon), 
$$
if $N\ge N_0$. Choose $(x_k^N)_{k=1}^{N-1}\subseteq\Bbb
A_N(s-\epsilon)$ as in the claim. Clearly,
$$
\align
&\frac 1{N^2}\log
Z_{N-1,N}^{\lambda,\sigma_N}(s_N)\ge-\frac{\beta}2\left[
\frac 1{N^2}\sum_{j\neq k}\log|x_j^N-x_k^N|^{-1}\right.\\&\left.+
\sum_{j=1}^{N-1} V_N(x_j^N)-\frac 1{N^2}\sum_{j=1}^{N-1}\log
|\sigma_N-x_j^N|\right],
\endalign
$$
and consequently, by Claim 4.2,
$$
\liminf_{N\to\infty}\frac 1{N^2}\log
Z_{N-1,N}^{\lambda,\sigma_N}(s_N)\ge-\frac{\beta}2 F_V^{s-\epsilon}.
$$
Since $F_V^{s-\epsilon}\searrow F_V^s$ as $\epsilon\to 0+$,
$$
\liminf_{N\to\infty}\frac 1{N^2}\log
Z_{N-1,N}^{\lambda,\sigma_N}(s_N)\ge-\frac{\beta}2 F_V^{s}.\tag 4.11
$$
Thus, given $\delta >0$, we can choose $N(\delta)$ so that if $N\ge
N(\delta)$ ,then
$$
\frac 1{N^2}\log
Z_{N-1,N}^{\lambda,\sigma_N}(s_N)\ge-\frac{\beta}2
(F_V^{s}+\delta).\tag 4.12
$$
It follows from (4.2) with $M=N-1$ and (4.3),that for any 
$0<\rho<1/2$,
$$
\align
&P_{N-1,N}^{\lambda,\sigma_N}(\Omega_{N,\eta+a}(S_N)^c;s_N)\\
&\le e^{\frac{\beta N^2}2 (F_V^s+\delta)}\sum_{x\in\Bbb
A(s_N)^{N-1}\setminus\Omega_{N,\eta+a}(s_N)}
e^{-\frac{\beta}2K_{N-1,V_N}(x) -\frac{\beta}2\sum_j
V_N(x_j)}\prod_{j=1}^{N-1}|\sigma_N-x_j|^{\lambda\beta}\\
&\le e^{\frac{\beta N^2}2(F_V^s+\delta)-\frac{\beta}2
(1-\rho)(F_V^s +\eta+a)N^2}\left[\sum_{t\in\Bbb
A_N}(t^2+1)^{-\frac{\beta}2\xi (N-1)}(1+\sigma_N^2)^{\lambda\beta/2}
\right]^N\\ 
&\le e^{-\frac{\beta}4 aN^2}
\endalign
$$
if $N$ is sufficiently large (independent of $a\ge 0$). Note that
$\delta+\rho F_V^s-\eta<0$ if we choose $\delta=\eta/2$ and $\rho$
sufficiently small. This completes the proof.
\enddemo

This lemma can be used to prove

\proclaim{Corollary 4.3} For any $s\in (1,\infty]$,
$$
\lim_{N\to\infty}\frac 1{N^2}\log Z_N(s)=-\frac {\beta}2F_V^s.\tag 4.13
$$
Furthermore $F_V^s-F_V>0$ if $s<b_V$. 
\endproclaim

\demo{Proof} The lower limit follows by taking
$\lambda=\sigma_N=0$ in (4.11) (replacing $N-1$ by $N$ does not modify
the argument above in any essential way). 
Given $0<\rho <1$, we can use
(4.2) with $M=N$ and the continuity of $\exp K_{N,V_N}$ to see that
$$
\align
Z_N(s)&=\sum_{x\in\Bbb
A_N(s)^N}e^{-\frac{\beta}2K_{N,V_N}(x)-\frac{\beta}2\sum_{j=1}^N
V_N(x_j)}\\
&\le\sup_{x\in\Bbb
A_N(s)^N}e^{-\frac{\beta}2(1-\rho)K_{N,V_N}(x)}\sum_{x\in\Bbb
A_N(s)^N}e^{-\frac{\beta}2\rho\xi(N-1)\sum_j\log(1+x_j^2)}\\
&\le e^{-\frac{\beta}2(1-\rho)K_{N,V_N}(y^N)+CN},\tag 4.14
\endalign
$$
if $N$ is sufficiently large, where $y^N=(y_1^N,\dots,y_N^N)\in\Bbb
A_N(s)^N$. Clearly, $y_j^N\neq y_k^N$ if $j\neq k$. Set
$\lambda_N=N^{-1}\sum_j\delta_{y_j^N}$. It follows from (4.12), with
$\lambda=\sigma=0$ and $N-1$ replaced by $N$, that $N^{-2}\log Z_N(s)\ge
-\beta (F_V^s+\delta)/2$ for $N\ge N(\delta)$, so (4.2) and (4.14)
yield
$$
\int_0^\infty\log(1+t^2)d\lambda_N(t)\le C.
$$
Thus $\{\lambda_N\}_{N=1}^\infty$ is tight. Pick a subsequence that
gives the upper limit of \linebreak $N^{-2}\log Z_N(s)$, and a further
subsequence so that $\lambda_{N_j}$ converges weakly to $\nu=\psi dx$.
The measyre $\nu$ has to be absolutely continuous with density
satisfying $0\le\psi\le 1$ because of the definition of $\lambda_N$.
Using (4.1) and $|t-s|\le\sqrt{t^2+1}\sqrt{s^2+1}$ we see that
$k_{V_N}(t,s)\ge 0$. Set, for given $M>0$,
$k_{V_N}^M(t,s)=\min(k_{V_N}(t,s),M)$
and choose $\phi_T(t)$ continuous so that $0\le\phi_T\le 1$,
$\phi_T(t)=1$ if $|t|\le T$, $=0$ if $|t|\ge T+1$ and $\phi_T(t)\le
\phi_{T'}(t)$ if $T\le T'$. Then, $k_{V_N}(t,s)\ge
\phi_T(t)\phi_T(s)k_{V_N}^M(t,s)$ and using the estimate (4.14) we get
$$\align
&\frac 1{N_j^2}\log Z_{N_j}(s)\\&\le\frac {C+\frac{\beta}2(1-\rho)M}{N_j}
-\frac{\beta}2(1-\rho)
\int_0^\infty\int_0^\infty\phi_T(t)\phi_T(s)K^M_{V_N}(t,s)
d\lambda_{N_j}(t)d\lambda_{N_j}(s),
\endalign$$
and thus, letting $j\to\infty$, $M\to\infty$, $T\to\infty$ and
$\rho\to 0+$ in that order, we obtain
$$
-\frac{\beta}2  F_V^s\le\liminf_{N\to\infty}\frac 1{N^2}\log Z_N(s)
\le\limsup_{N\to\infty}\frac 1{N^2}\log Z_N(s)\le
-\frac{\beta}2I_V[\psi].
$$
Thus $I_V[\psi]\le F_V^s$ and $\psi\in\Cal A_s$, 
so we must have $\psi=\phi_V^s$.

Assume that $F_V^s\le F_V$ and $s<b_V$. Then $I_V[\phi_V^s]\le
I_V[\phi_V]$ and consequently $\phi_V^s=\phi_V$ by the uniqueness of
the minimizing measure. This contradicts the definition of $b_V$. The
corollary is proved.
\enddemo

Note that by (2.13) corollary 4.3 implies (2.16) so we have proved the
first part of theorem 2.2. Before turning to the proof of the second
part we need one more consequence of lemma 4.1.

\proclaim{Corollary 4.4} Let $\{s_N\}$ be as in lemma 4.1 and assume
that $f:[0,\sigma+\epsilon]\to\Bbb R$, $\epsilon >0$, is continuous,
or $f:[0,\infty)\to\Bbb R$ is continuous and bounded in case
$s_N\equiv\infty$. Then
$$
\lim_{N\to\infty}\frac 1N\log E_{N-1,N}^{y,\sigma_N}[e^{\sum_{j=1}^N
f(x_j)};s_N] =\int_0^\infty f(t)\phi_V^\sigma(t)dt.\tag 4.15
$$
Furthermore let
$$
u_{N-1,N}^{y,\sigma_N}(t)=\frac 1{N-1}
E_{N-1,N}^{y,\sigma_N}[\sum_{i=1}^{N-1}\delta_{t,x_i}],\tag 4.16
$$
($\delta_{t,s}$ is Kronecker's delta), be the 1-dimensional marginal
distribution of the probability measure (4.4) (with $s_N\equiv\infty$).
Then for each $0<y\le 1$:
\item{(i)} $0\le u_{N-1,N}^{y,\sigma_N}(t)\le\frac 1{N-1}$
for all $t\in\Bbb A_N$,
\item{(ii)} if $\delta_t$ is the Dirac measure at $t$, then
$\sum_{t\in\Bbb A_N}u_{N-1,N}^{y,\sigma_N}(t)\delta_t$
converges weakly to $\phi_V(t)dt$ as $N\to\infty$.
\item{(iii)} $u_{N-1,N}^{y,\sigma_N}(\sigma_N)=0$.
\endproclaim

\demo{Proof} We can prove (4.15) using lemma 4.1 in exactly the same
way as in the proof of (2.5) on p. 194 in [Jo], see also [De]. The
weak limit (ii) is a direct consequence of (4.15), see [De]. Note that
the limit does not depend on $y$ since the factor
$\prod_{i=1}^{N-1}|\sigma_N-x_i|^{y\beta}$ does not affect the leading
asymptotics. 

In the expectation (4.16) all the $x_i\,$:s have to be different,
otherwise the probability is zero, and consequently the expectation is
$\le 1$, which proves (i). That (iii) holds follows from the presence
of the factor $\prod_{i=1}^{N-1}|\sigma_N-x_i|^{y\beta}$. The
corollary is
proved.
\enddemo

We turn now to the proof of the 
upper-tail limit. Note that
$$
Q_{M,N}(x)=e^{-\frac{N\beta}2V_N(x_M)}\prod_{i=1}^{M-1}|x_M-x_i|^\beta
Q_{M-1,N}(x'),\tag 4.17
$$
where $x'=(x_1,\dots,x_{M-1})$. Using this identity we see that
$$
\align
Z_{M,N}(t)&=M!\sum\Sb x\in\Bbb A_N^M\\x_1\le \dots\le x_M\le t\endSb
Q_{M,N}(x)\\
&=M\sum_{s\in\Bbb A_N(t)}e^{-\frac{N\beta}2V_N(s)}
\sum_{x\in\Bbb A_N(s)^{M-1}}\prod_{i=1}^{M-1}|s-x_i|^\beta
Q_{M-1,N}(x).
\endalign
$$
If we define
$$
H_{M-1,N}(s)=\frac 1{Z_{M-1,N}(s)}
\sum_{x\in\Bbb A_N(s)^{M-1}}\prod_{i=1}^{M-1}|s-x_i|^\beta
Q_{M-1,N}(x)
$$
this can be written
$$
Z_{M,N}(t)=M\sum_{s\in\Bbb A_N(t)}e^{-\frac{N\beta}2V_N(s)}
Z_{M-1,N}(s)H_{M-1,N}(s),\tag 4.18
$$
or
$$
F_{M,N}(t)=\frac{MZ_{M-1,N}}{Z_{M,N}}
\sum_{s\in\Bbb A_N(t)}e^{-\frac{N\beta}2V_N(s)}
F_{M-1,N}(s)H_{M-1,N}(s).\tag 4.19
$$
This is the main formula to be used in the proof of (2.17). We will
need two choices of $M$, namely $M=N$ and $M=N-1$. They are handled
completely analogously and we will consider only the case $M=N$.

Write $\Bbb A_N(t,s)=\Bbb A_N\cap(t,s)$ for any $0\le t<s\le \infty$
and $\Bbb A_N(t)^\ast=\Bbb A_N(t,\infty)$.
If we let $t\to\infty$ in (4.19) and then subtract (4.19) from the
limiting equality, we get
$$
1-F_N(t)=\frac{NZ_{N-1,N}}{Z_{N,N}}
\sum_{s\in\Bbb A_N(t)^\ast}e^{-\frac{N\beta}2V_N(s)}
F_{N-1,N}(s)H_{N-1,N}(s).\tag 4.20
$$
Set
$$
\Phi_V=F_V-\frac 12\int_0^\infty V(s)\phi_V(s)ds,
$$
From the variational relations for $\phi_V(t)$ it follows that
$$
\int_0^\infty\log|b_V-s|^{-1}\phi_V(s)ds+\frac 12 V(b_V)=\Phi_V.\tag
4.21
$$
\proclaim{Lemma 4.5} We have
$$
\limsup_{N\to\infty}\frac 1N\log\frac{Z_{N-1,N}}{Z_{N,N}}\le\beta\Phi_V.\tag
4.22
$$
\endproclaim
\demo{Proof}
By (4.17) we have
$$
\align
\frac{Z_{N,N}}{Z_{N-1,N}}&=
\sum_{s\in\Bbb
A_N}e^{-\frac{N\beta}2V_N(s)}E_{N-1,N}[\prod_{i=1}^{N-1}
|s-x_i|^\beta]\\
&\ge e^{-\frac{N\beta}2V_N(r)}E_{N-1,N}[\prod_{i=1}^{N-1}
|r-x_i|^\beta]\tag 4.23
\endalign
$$
for any $r\in\Bbb A_N$. One difficulty in estimating the right hand
side in (4.23)
comes from the fact that, due to the discrete nature of the problem
the integrand could, apriori, be zero for many $y\,$:s with high
probability. Note that we define $0^y=0$ for any $y>0$.
Let $\psi_s(t)=1$ if $t\neq s$ and $\psi_s(s)=0$. 

Consider
$$
f_N(y;s)=\frac 1N\log E_{N-1,N}[\prod_{i=1}^{N-1}
|s-x_i|^{y\beta}\psi_s(x_i)].
$$
Then,
$$\align
f_N(0+;s)&=\lim_{y\to 0+} f_N(y;s)=\frac
1N\log E_{N-1,N}[\prod_{i=1}^{N-1}\psi_s(x_i)] \\
&=\frac 1N\log P_{N-1,N}[\text{all $x_i\neq s$}].\tag 4.24
\endalign
$$
Let $\epsilon >0$ be given and write $\Bbb B_N(\epsilon)=\Bbb
A_N(b_V+\epsilon,b_V+2\epsilon)$. Now,
$$\align
\sum_{s\in\Bbb B_N(\epsilon)}P_{N-1,N}[\text{all $x_i\neq s$}]&\ge
P_{N-1,N}[\bigcup_{s\in\Bbb B_N(\epsilon)}\{\text{all $x_i\neq
s$}\}]\\
&=1-P_{N-1,N}[\bigcap_{s\in\Bbb B_N(\epsilon)}\{\text{one $x_i=
s$}\}].\tag 4.25
\endalign
$$
Take $g:[0,\infty)\to[0,\infty)$ continuous such that $g(s)=1$ if 
$b_V+\epsilon\le s\le b_V+2\epsilon$ and $g(s)=0$ if $0\le s\le b_V$ or
$s\ge b_V+3\epsilon$. Then,
$$
e^{\epsilon N}P_{N-1,N}[\bigcap_{s\in\Bbb B_N(\epsilon)}\{\text{one $x_i=
s$}\}]\le E_{N-1,N}[e^{\sum_{i=1}^N g(x_i)}]\le e^{\epsilon N/2}\tag
4.26
$$
for all sufficiently large $N$. The first inequality follows from the
definitions whereas the second follows from corollary 4.4,
(4.15). Combining (4.25) and (4.26) we see that
$$
\max_{s\in\Bbb B_N(\epsilon)} P_{N-1,N}[\text{all $x_i\neq s$}]\ge
\frac 1{2N}\tag 4.27
$$
for all sufficiently large $N$. Hence, by (4.24) and (4.27) we can
choose $\sigma_N=\sigma_N(\epsilon)\in\Bbb B_N(\epsilon)$ so that
$$
\lim_{N\to\infty}f_N(0+;\sigma_N)=0.\tag 4.28
$$

Take $r=\sigma_N$ in (4.23). Then
$$\align
\frac 1N\log\frac{Z_{N,N}}{Z_{N-1,N}}&\ge
-\frac{\beta}{2}V_N(\sigma_N)+
f_N(1;\sigma_N)\\
&=-\frac{\beta}{2}V_N(\sigma_N)+
f_N(0+;\sigma_N)+\beta\int_0^1 f_N'(y;\sigma_N)dy.\tag 4.29
\endalign$$
We can pick a subsequence $\{N_j\}$ which gives $\liminf_{N\to\infty}
\frac 1N\log\frac{Z_{N,N}}{Z_{N-1,N}}$ and such that
$\sigma_{N_j}(\epsilon)\to\sigma(\epsilon)\in
[b_V+\epsilon,b_V+2\epsilon]$. Then, by (4.28) and (4.29),
$$
\liminf_{N\to\infty}
\frac 1N\log\frac{Z_{N,N}}{Z_{N-1,N}}\ge -\frac {\beta}2
V(\sigma(\epsilon)) +\beta\liminf_{j\to\infty}
\int_0^1 f_{N_j}'(y;\sigma_{N_j})dy.\tag 4.30
$$
Now,
$$\align
f_N'(y;\sigma_N)&=E_{N-1,N}^{y,\sigma_N}[\frac 1N\sum_{i=1}^{N-1}
\log|\sigma_N-x_i|] \\
&=\frac{N-1}N\sum_{t\in\Bbb
A_N}\log|\sigma_N-t|u_{N-1,N}^{y,\sigma_N}(t).
\endalign
$$
Hence, by corollary 4.4 (i) and (iii),
$$
f_N'(y;\sigma_N)\ge 2\frac{1}N\sum_{i=1}^N\log\frac iN\ge -2
$$
and consequently, by Fatou's lemma,
$$
\liminf_{j\to\infty}\int_0^1 f_{N_j}'(y;\sigma_{N_j})dy\ge
\int_0^1 \liminf_{j\to\infty}f_{N_j}'(y;\sigma_{N_j})dy.\tag 4.31
$$

Given $\delta>0$, small, and  $M>0$ set
$$
f_{M,\delta}(t)=\cases \log M, & \text{if $|t|\ge M$}\\
\log |t|, & \text{if $\delta\le |t|< M$}\\
\log \delta, & \text{if $|t|\le\delta$}.
\endcases
$$
By corollary 4.4  (i) and (iii) we have
$$
\align
&\left|\sum_{t\in\Bbb A_N}(\min(\log
M,\log|\sigma_N-t|)-f_{M,\delta}(\sigma_N
-t))
u_{N-1,N}^{y,\sigma_N}(t)\right|\\
&\le\sum_{t\in\Bbb
A_N\,;\,0<|t-\sigma_N|\le\delta}\biggl|\log\biggl|\frac{\sigma_N-t}{\delta} 
\biggr|\biggr|\frac 1{N-1}\le\frac 2{N-1}\sum_{k=1}^{[N\delta]}
\log\frac{N\delta}k \\
&\le\frac{2N}{N-1}\delta.
\endalign
$$
Also, if $|\sigma_N-\sigma_\epsilon|\le\delta$, which is true if $N$
is large enough
$$
|f_{M,\delta}(|\sigma_N-t|)-f_{M,\delta}(|\sigma(\epsilon)-t|)|
\le\delta\log\frac1{\delta} .
$$
Since $\log|\sigma_N-t|\ge\min (\log M,\log|\sigma_N-t|)$ and $M$,$\delta$ are
arbitrary it follows from corollary 4.4, (ii) that
$$
\liminf_{j\to\infty}f_{N_j}'(y;\sigma_{N_j})\ge\int_0^\infty
\log|\sigma(\epsilon)-t|\phi_V(t)dt.
$$
Together with (4.30) and (4.31) this gives
$$
\liminf_{N\to\infty}\frac 1N\log\frac{Z_{N,N}}{Z_{N-1,N}}\ge
-\frac{\beta}2 V(\sigma(\epsilon))+\beta\
\int_0^\infty\log|\sigma(\epsilon)-t|\phi_V(t)dt.
$$
We can pick a sequence $\epsilon_j\to 0$ such that
$\sigma(\epsilon_j)\to b_V$ and using (4.24) we obtain
$$
\liminf_{N\to\infty}\frac 1N\log\frac{Z_{N,N}}{Z_{N-1,N}}\ge
-\beta\Phi_V,
$$
and the lemma is proved.
\enddemo

Given $\delta>0$ we can use lemma 4.5 to find $N_0(\delta)$ so that 
$$
\frac{Z_{N-1,N}}{Z_{N,N}}\le e^{N\beta(\Phi_V+\delta)}\tag 4.32
$$
if $N\ge N_0(\delta)$. Since $F_{N-1,N}(s)\le 1$ we can combine (4.20)
and (4.32) to get the estimate
$$
1-F_N(t)\le Ne^{N\beta(\Phi_V+\delta)}\sum_{s\in\Bbb A_N(t)^\ast}
e^{-\frac {N\beta}2 V_N(s)}H_{N-1,N}(s).\tag 4.33
$$
We have
$$
\align
H_{N-1,N}(s)&=E_{N-1,N}[\prod_{i=1}^{N-1}|s-x_i|^\beta;s]\\
&\le(1+s^2)^{\frac{\beta}2(N-1)}
E_{N-1,N}^{0,0}[\prod_{i=1}^{N-1}(1+x_i^2)^{\beta/2};s]\le e^{CN}
(1+s^2)^{\beta N/2},
\endalign
$$
where the last inequality is proved, using lemma 4.1, just as (4.25)
in [Jo]. Together with (4.1) this gives
$$
e^{-\frac {N\beta}2 V_N(s)}H_{N-1,N}(s)\le e^{CN-\frac{N\beta\xi}2\log
(1+s^2)},\tag 4.34
$$
Hence, given a constant $D>0$, there is a constant $d>0$ such that
$$
e^{N\beta(\Phi_V+\delta)}\sum_{s\in\Bbb A_N(d)^\ast}e^{-N\beta
V_N(s)/2} H_{N-1,N}(s)\le e^{-ND}.\tag 4.35
$$

For $t\ge s$ we define
$$
H_{N-1,N}(t,s)=\frac 1{Z_{N-1,N}(s)}\sum_{x\in\Bbb A_N(s)^{N-1}}
\prod_{j=1}^{N-1}|t-x_i|^\beta Q_{N-1,N}(x).
$$
Clearly,
$$
H_{N-1,N}(s)= H_{N-1,N}(s,s)\le H_{N-1,N}(t,s)\tag 4.36
$$
if $t\ge s$.
Combining the estimates (4.33), (4.35) and (4.36) we obtain
$$
1-F_N(t)\le Ne^{-ND}+Ne^{N\beta(\Phi_V+\delta)}
\sum_{x\in\Bbb A_N(t,d)}e^{-\frac {N\beta}2
V_N(s)}H_{N-1,N}(s+\epsilon,s)
\tag 4.37
$$
for any $\epsilon>0$. Let $s_N\in
\Bbb A_N(t,d)$ be the $s$ which gives the largest term in the sum in
(4.37). Then
$$
1-F_N(t)\le Ne^{-ND}+N^2(d-t)e^{N\beta(\Phi_V+\delta-\frac
12V_N(s_N))}H_{N-1,N}(s_N+\epsilon,s_N).\tag 4.38
$$
Choose a sequence which gives the upper limit of $N^{-1}\log
(1-F_N(t))$ and such that $s_{N_j}\to\sigma\in [t,d]$. We would like
to prove that
$$
\lim_{j\to\infty}\frac 1{N_j}\log H_{N_j-1,N_j}(s_{N_j}+\epsilon,s_{N_j})=
-\beta\int\log|\sigma+\epsilon-t|\phi_V^\sigma(t)dt.\tag 4.39
$$
We will write $N$ instead of $N_j$ for simplicity.
Looking at the definition of \linebreak $H_{N-1,N}(t,s)$, we see that we are
interested in the limit of
$$
\frac 1N\log E_{N-1,N}[e^{\beta\sum_{j=1}^{N-1}\log
|s_N+\epsilon-x_i|} ;s_N]
$$
as $N\to\infty$, $s_N\to\sigma$. Since
$$
|\log |s_N+\epsilon-x_i|-\log |\sigma+\epsilon-x_i||
=|\log|1+\frac{s_N-\sigma}{\sigma+\epsilon-x_i}||\le C
\frac{|s_N-\sigma|}{\sigma+\epsilon-x_i},\tag 4.40
$$
where $C$ is a numerical constant, and $s_N\le\sigma+\epsilon/2$ for
$N$ large enough, the limit (4.39) follows from corollary 4.4. 

If $t>b_V$, then $\phi_V^\sigma=\phi_V$, since $\sigma\ge t$, and
combining (4.38) and (4.39) yields
$$\align
&\limsup_{N\to\infty}\frac 1N\log(1-F_N(t))\\
&\le\max\{-D,\beta\Phi_V+\delta-\frac{\beta}2 V(\sigma)-\beta\int\log
|\sigma+\epsilon-t|^{-1}\phi_V(t)dt\}.\tag 4.41
\endalign
$$
Note that $\sigma$ could depend on $\epsilon$ and $d$. Pick a sequence
$\epsilon=\epsilon_j\to 0+$ and then a subsequence so that
$\sigma(\epsilon_{j_k}) \to\tau\in [t,d]$. Then, since $D$ and
$\delta$ are arbitrary, we get
$$
\limsup_{N\to\infty}\frac 1N\log(1-F_N(t))\le \beta(\Phi_V-\inf_{\tau\ge t}\int
k_V(\tau,s)\phi_V(s)ds)\tag 4.42
$$
and we have proved one half of (2.17).

We now turn to the lower limit. If we start with $M=N-1$ instead of
$N$ then (4.42) holds with $F_{N-1}$ replaced by $F_{N-1,N}(t)$. By
assumption the right hand side of (4.42) is negative for all
$t>b_V$. Hence, if $t>b_V$, we see that
$$
F_{N-1,N}(t)\ge 1/2\tag 4.43
$$
for all sufficiently large $N$. Note that, if $t\ge s$, then
$$
H_{N-1,N}(t)\ge \frac{Z_{N-1,N}(s)}{Z_{N-1,N}(t)}H_{N-1,N}(t,s)
\ge F_{N-1,N}(s)H_{N-1,N}(t,s).\tag 4.44
$$

The function $f(\tau)=\int k_V(\tau,s)\phi_V(s)ds$ is continuous 
on $[t,\infty)$ and
$f(\tau)\to\infty$ as $\tau\to\infty$, so it assumes its minimum in
$[t,\infty)$ at some point $\tau_0\ge t$. Let $\epsilon>0$. Pick
$s_N\in \Bbb A_N(t)^\ast$ such that
$s_N\searrow\tau_0+\epsilon$. Then, picking one term in the sum
$$\align
&\sum_{s\in \Bbb A_N(t)^\ast}e^{-\frac{N\beta}2 V_N(s)}
F_{N-1,N}(s)H_{N-1,N}(s)\\
&\ge e^{-\frac{N\beta}2 V_N(s_N)}
F_{N-1,N}(\tau_0)^2H_{N-1,N}(s_N,s_N-\epsilon).
\endalign$$
If we use the limit (4.39), the estimate (4.43) with $s=\tau_0$, and
let $\epsilon\to 0+$, we see that
$$
\align
&\liminf_{N\to\infty}\frac 1N\log
\sum_{s\in \Bbb A_N(t)^\ast}e^{-\frac{N\beta}2 V_N(s)}
F_{N-1,N}(s)H_{N-1,N}(s)\\
&\ge-\frac{\beta}2
V(\tau_0)-\beta\int\log|\tau_0-t|^{-1}\phi_V(t)dt.\tag 4.45
\endalign
$$
To complete the proof we need

\proclaim{Lemma 4.6} For any $V_N$ satisfying the conditions (i) -
(iii),
$$
\liminf_{N\to\infty}\frac 1N
\log\frac{Z_{N-1,N}}{Z_{N,N}}\ge\beta\Phi_V.\tag 4.46
$$
\demo{Proof} If we let $t\to\infty$ in (4.19), we see that,
$\epsilon>0$,
$$\align
\frac{Z_{N,N}}{Z_{N-1,N}}&=N\sum_{s\in\Bbb A_N}e^{-\frac{N\beta}2 V_N(s)}
F_{N-1,N}(s)H_{N-1,N}(s)\\
&\le N\sum_{s\in\Bbb A_N(b_V-\epsilon)}e^{-\frac{N\beta}2 V_N(s)}
F_{N-1,N}(s)H_{N-1,N}(s)\\
&+N\sum_{s\in\Bbb A_N(b_V-\epsilon)^\ast}e^{-\frac{N\beta}2 V_N(s)}
H_{N-1,N}(s),\tag 4.47
\endalign
$$
since $F_{N-1,N}(s)\le 1$. By adjusting the constant $C$ we see that
(4.34) holds for all $s\in\Bbb A_N$, so the first sum in the right
hand side of (4.47) is
$$
\le e^{CN}F_{N-1,N}(b_V-\epsilon)\sum_{s\in\Bbb
A_N}e^{-\frac{\beta}2N\xi\log(1+s^2)} \le
e^{CN-\frac{\beta}2L(b_V-\epsilon)N^2} 
$$
for all sufficiently large $N$ by the first part of theorem 2.2.
(Replacing $F_N(t)$ by $F_{N-1,N}(t)$ does not make any difference.)
Since $L(b_V-\epsilon)>0$ if $\epsilon>0$, the first part of the right
hand side of (4.47) is negligible.

The same argument that lead us from (4.33) to (4.42) allows us to
treat the second term in the right hand side of (4.47) and obtain
$$
\align
&\limsup_{N\to\infty}\frac 1N\log\frac{Z_{N,N}}{Z_{N-1,N}}\\
&\le\max\{-D,-\frac{\beta}2 V(\sigma)-\beta\int\log
|\sigma+\eta-t|^{-1}\phi_V^\sigma(t)dt\}.\tag 4.48
\endalign
$$
where $\sigma\in [b_V-\epsilon,d]$, $\eta>2\epsilon$, $D>0$ are given.
Take $\epsilon=\epsilon_j\to 0+$ so that $\sigma(\epsilon_j)\to\tau\in
[b_V,d]$ . Note that $\phi_V^{\sigma(\epsilon_j)}(t)dt$ converges
weakly to $\phi_V^\tau(t)dt=\phi_V(t)dt$. Using an inequality like
(4.40) we get
$$
\align
&\limsup_{N\to\infty}\frac 1N\log\frac{Z_{N,N}}{Z_{N-1,N}}\\
&\le\max\{-D,-\frac{\beta}2 V(\tau)-\beta\int\log
|\tau+\eta-t|^{-1}\phi_V(t)dt\}.\tag 4.49
\endalign
$$
We can now repeat the argument that lead from (4.41) to (4.42) and
obtain
$$
\limsup_{N\to\infty}\frac 1N\log\frac{Z_{N,N}}{Z_{N-1,N}}
\le \frac{\beta}2\int V(s)\phi_V(s)ds-\beta\inf_{\tau\ge b_V}\int k_V(\tau,
s)\phi_V(s)ds\le -\beta\Phi_V,
$$
since $\int k_V(\tau,s)\phi_V(s)ds\ge F_V$ if $\tau\ge b_V$. The lemma
is proved.
\enddemo

Combining (4.20), (4.45) and lemma 4.6, we see that
$$\align
&\liminf_{N\to\infty}\frac 1N\log(1-F_N(t))\\
&\ge\beta (F_V-\int k_V(\tau_0,s)\phi_V(s)ds)=
\beta (F_V-\inf_{\tau \ge t}\int k_V(\tau,s)\phi_V(s)ds),
\endalign
$$
by the choice of $\tau_0$. This completes the proof of theorem 2.2.
\enddemo
\bigskip
\heading 5. Asymptotics for the Meixner kernel\endheading
\medskip
This section is devoted to the proof of lemma 3.2, which is based on
establishing the appropriate asymptotics of the Meixner polynomials.
See [Go] and [JW] for some results on the 
asymptotics of Meixner polynomials. 

From (3.3) we obtain, $x\in\Bbb R$,
$$
\align
m_n^{K,q}(x)&=(-1)^n\frac{(\sqrt{\gamma})^{n+K}n!}{(\sqrt{q})^{n}2\pi i}
\int_{\Gamma_r}\biggl(\frac{\sqrt{\gamma}+z/\sqrt{q}}{\sqrt{\gamma}+\sqrt{q} 
z}\biggr)^x\frac{dz}{(\sqrt{\gamma}+\sqrt{q}z)^Kz^{n+1}}\\
&-\frac{\sin\pi x}{\pi} 
\frac{(\sqrt{\gamma})^{n+K}n!}{(\sqrt{q})^{n}}\int_{\sqrt{\gamma q}}^r
\biggl|\frac{\sqrt{\gamma}-t/\sqrt{q}}{\sqrt{\gamma}-\sqrt{q} 
t}\biggr|^x\frac{dt}{(\sqrt{\gamma}-\sqrt{q}t)^Kt^{n+1}},\tag 5.1
\endalign
$$
where $\Gamma_r$ is the circle $|z|=r$, $0<r<\sqrt{\gamma/q}$; if
$0<r\le\sqrt{\gamma q}$ the second integral should be omitted. Let
$b=(1+\sqrt{\gamma q})^2/(1-q)$ as before, let $\sigma$ be given by
(1.11) and set
$$
a=b+\gamma-1=\frac{(\sqrt{\gamma}+\sqrt{q})^2}{1-q}.
$$
Set
$$\align
t(z)&=\biggl(\frac{\sqrt{\gamma q}+z}{\sqrt{\gamma q}+1}\biggr)
\biggl(\frac{\sqrt{\gamma}+\sqrt{q}}{\sqrt{\gamma}+\sqrt{q}z}\biggr),\\
s(z)&=\frac{\sqrt{\gamma}+\sqrt{q}}{\sqrt{\gamma}+\sqrt{q}z},
\endalign
$$
and
$$
A_N(x)=\frac{b^x}{x^{x+K}}\frac{(x+K-1)!N!}{x!(N+K-2)!}
\frac{\gamma^{K+N}}{1-q} \sqrt{\frac{q}{\gamma}}.
$$
For $0<r<\sqrt{\gamma/q}$ we define
$$
D_n^r(x;g)=\frac 1{2\pi}\int_{-\pi}^\pi
g(re^{i\theta})t(re^{i\theta})^xs(re^{i\theta})^K
\frac{d\theta}{r^ne^{in\theta}},\tag 5.3
$$
$F_n^r(x;g)=0$ if $r\ge\sqrt{\gamma q}$ and
$$
F_n^r(x;g)=(-1)^{n+x+1}\int_{\sqrt{\gamma
q}}^r|t(-\tau)|^xs(-\tau)^Kg(-\tau) \frac{d\tau}{\tau^{n+1}}.\tag 5.4
$$
The powers are defined by taking the prinipal branch of the logarithm.

The Meixner kernel (3.4) can now be written, for $x,y$ {\it integers }
(which is the case we need),
$$
K_N(x,y)=\sqrt{A_N(x)A_N(y)}\frac{D_N(x;g_1)D_N(y;g_2)-D_N(x;g_2)D_N(y;g_1)} 
{x-y}\tag 5.5
$$
if $x\neq y$, and
$$\align
K_N(x,x)&=A_N(x)[D_N(x-1;g_3)D_N(x;g_2)- D_N(x;g_1)D_N(x-1;g_4)\\
&+F_N(x;g_1)D_N(x;g_2)-F_N(x;g_2)D_N(x;g_1)],\tag 5.6
\endalign$$
where $g_1(z)\equiv 1$, $g_2(z)=z-1$, $g_3(z)=t(z)\log t(z)$ and
$g_4(z)=g_2(z)g_3(z)$. The functions $g_i(z)$ are bounded for $|z|\le
1$.

Write $x=Nb+y$ and $K=[\gamma N]-N+1\doteq N(\gamma_N-1)\doteq
N(\gamma-1)+\omega_N$, $0<\omega_N\le 1$.

\proclaim{Lemma 5.1} If $x=Nb+\xi\sigma N^{1/3}$ and $M_0>0$ is a
given constant, there are constants $c_1(q,\gamma)$ and
$c_2(q,\gamma)$, such that
$$
\frac 1N A_N(x)\le c_1(q,\gamma)e^{c_2(q,\gamma)\xi N^{-2/3}}\tag 5.7
$$
for all $\xi\ge -M_0$. Furthermore,
$$
\lim_{N\to\infty}\frac 1N
A_N(x)=\frac{\gamma\sqrt{q}}{(1-q)\sqrt{ab}}\tag 5.8
$$
uniformly for $|\xi|\le M_0$.
\endproclaim

\demo{Proof} By Stirling's formula
$$\align
A_N(x)&=\frac{(x+K)^{x+K}N^Nb^x}{x^x(N+K)^{N+K}a^{x+K}}\gamma
^{K+N}\frac {(N+K)(N+K-1)}{x+K}\\
&\times\sqrt{\frac{(x+K)N}{x(N+K)}}\frac 1{1-q}\sqrt{\frac
q{\gamma}}e^{o(1)}.\tag 5.9
\endalign
$$
Write $a_N=b+\gamma_N-1$. Then,
$$
\frac{(x+K)^{x+K}N^Nb^x}{x^x(N+K)^{N+K}a^{x+K}}\gamma
^{K+N}=\biggl(\frac{Nb}x\biggr)^x\biggl(\frac{x+K}{Na_N}\biggr)^{x+K}
\biggl(\frac{a_N}a\biggr)^{x+K}\biggl(\frac{\gamma}{\gamma_N}\biggr)^{N+K}.\tag
5.10
$$
If we write
$u=Na_N$ and $v=Nb<u$. Then
$$
\biggl(\frac{Nb}{x}\biggr)^{x}
\biggl(\frac{x+K}{Na_N}\biggr)^{x+K}=
\biggl(1+\frac yu\biggr)^{u+y}\biggl(1+\frac yv\biggr)^{-v-y}\doteq
e^{g(y)}. 
$$
Since $g(0)=g'(0)=0$ and $g''(t)=(v-u)(u+t)^{-1}(v+t)^{-1}<0$, we have
$\exp{g(t)}\le 1$ if $\xi\ge 0$. If $-M_0\le\xi\le M_0$, then
$$
|g(t)|=|\int_0^t(t-s)g''(s)ds|\le CN^{-1/3}.
$$
Furthermore
$$
\biggl(\frac{a_N}{a}\biggr)^{x+K}=e^{\omega_N+O(\xi N^{-2/3})+o(1)}
$$
and
$$
\biggl(\frac{\gamma}{\gamma_N}\biggr)^{K+N}=e^{-\omega_N+o(1)}.
$$
Inserting these estimates into (5.10) we obtain
$$
\frac{(x+K)^{x+K}N^Nb^x}{x^x(N+K)^{N+K}a^{x+K}}\gamma
^{K+N}\le Ce^{C\xi N^{-2/3}}
$$
for $\xi\ge -M_0$ and
$$
\lim_{N\to\infty}\frac{(x+K)^{x+K}N^Nb^x}{x^x(N+K)^{N+K}a^{x+K}}\gamma
^{K+N}=1
$$
uniformly for $|\xi|\le M_0$. By (5.9) this proves (5.7) and
(5.8). The lemma is proved.
\enddemo

Set 
$$
u(z)=b\log(\sqrt{\gamma q}+z)-a\log(\sqrt{\gamma}+\sqrt{q}z)-\log z
$$
so that
$$
D_N^r(x;g)=\frac 1{2\pi}\int_{\Gamma_r}e^{N(u(z)-u(1))+y\log
t(z)+\omega_N\log s(z)}g(z)\frac {dz}{iz}.\tag 5.11
$$
Now,
$$\align
u'(z)&=-\rho(1-z)^2\\
&+\rho(1-z)^3\frac{\sqrt{q}z^2+(\sqrt{q}+\sqrt{\gamma}+
q\sqrt{\gamma})z+\sqrt{q}+
\sqrt{\gamma}+q\sqrt{\gamma}+\gamma\sqrt{q}} 
{z(z+\sqrt{\gamma q})(\sqrt{\gamma}+\sqrt{q}z)},
\endalign$$
where
$$
\rho=\frac{\gamma\sqrt{q}}{(1+\sqrt{\gamma
q})(\sqrt{\gamma}+\sqrt{q})}.
$$
Hence we can write
$$
u(z)-u(1)=\frac 13\rho(1-z)^3+\rho(1-z)^4v(z),\tag 5.12
$$
where one verifies that $|v(z)|\le 28/27$ if $|z-1|\le 1/4$.

By taking absolute values in (5.3) we obtain
$$
|D_N^r(x;g)|\le\frac {C}{2\pi}\biggl(\frac
ab\biggr)^{x/2}\frac{a^K(1-q)^K}{r^N}
\int_{-\pi}^{\pi}e^{f(\cos\theta)}d\theta,\tag 5.13
$$
where
$$
f(\tau)=\frac x2\log(\gamma q+r^2+2\sqrt{\gamma q}r\tau)+
\frac{x-K}2\log(\gamma +qr^2+2\sqrt{\gamma q}r\tau).
$$
Write $r=1-\delta$, $0\le\delta <1$. A computation shows that
$f'(\tau)\ge 0$ if (say)
$$
y\ge -\delta\frac{1+q+2\sqrt{\gamma q}}{1-q}N,\tag 5.14
$$
which covers all the $y\,$:s we are interested in.
Thus, if (5.14) is fullfilled, then
$$
|D_N^r(x;g)|\le C\exp(N(u(1-\delta)-u(1))+y\log t(1-\delta)).\tag 5.15
$$
By (5.12),
$$
u(1-\delta)-u(1)\le\rho\delta^3(\frac 13\delta\frac{28}{27})\le\frac
23\rho\delta^3\tag 5.16
$$
if $0\le\delta\le 1/4$. Now,
$$
\align
\log t(1-\delta)&=\log\bigl(1-\frac 1{1-\frac{\sqrt{q}}{\sqrt{\gamma}+
\sqrt{q}}\delta}\frac{(1-q)\sqrt{\gamma}}{(1+\sqrt{\gamma
q})(\sqrt{\gamma}+\sqrt{q})}\delta\bigl)\\
&\le-\rho(1-q)\frac 1{\sqrt{\gamma q}}\delta,
\endalign
$$
and consequently it follows from (5.15) and (5.16) that, if $y\ge 0$,
then
$$
|D_N^r(x;g)|\le C[\exp\bigl[\frac{2N}3\rho\delta^3-\rho
(1-q)\frac 1{\sqrt{\gamma q}}\delta y\bigr].\tag 5.17
$$
Recall that $y=\sigma N^{1/3}\xi$ with $\sigma$ given by (1.11).
Note that $\sigma=(1-q)^{-1}\sqrt{\gamma q}\rho^{-2/3}$. Choose
$\delta= (\rho N)^{-1/3}\sqrt{\xi}$ if $\xi\le (N\rho)^{2/3}/16$ and
$\delta=1/4$, if $\xi\ge (N\rho)^{2/3}/16$. Inserting this into (5.17)
gives
$$
|D_N^r(x;g)|\le C\exp\bigl[-\frac 13\min(\sqrt{\xi},\frac
14(N\rho)^{1/3})\xi \bigr],\tag 5.18
$$
for $\xi\ge 0$.

Let $\epsilon\in [0,\pi]$ and set
$$\align
I_1'&=\frac 1{2\pi}\int_{-\epsilon}^\epsilon g(re^{i\theta})
t(re^{i\theta})^xs(re^{i\theta})^K\frac{d\theta}{r^Ne^{iN\theta}},\\
I_1''&=D_N^r(x;g)-I_1'.
\endalign
$$
By the same argument that was used for (5.13) above, we see that if
$y$ satisfies (5.14), then
$$\align
|I_1''|&\le C|t(re^{i\epsilon})|^x|s(re^{i\epsilon})|^K\frac 1{r^N}\\
&\le C\exp\bigl[N\re
(u(re^{i\epsilon})-u(1))+y\log|t(re^{i\epsilon})|\bigr].
\tag 5.19
\endalign
$$

Next, we consider $F_N^r(x;g)$, $\sqrt{\gamma q}< r\le 1$. Taking
absolute values in (5.4) yields
$$
|F_N^r(x;g)|\le C\int_{\sqrt{\gamma q}}\biggl|\frac{\sqrt{\gamma
q}-\tau}{\sqrt{\gamma q}+1}\biggr|^x\biggl|\frac{\sqrt{\gamma}+\sqrt{q}}
{\sqrt{\gamma}-\sqrt{q}\tau}\biggr|^{x+K}\frac{d\tau}{\tau^{N+1}}.\tag
5.20
$$
The integrand in (5.20) is a increasing function of $\tau$ for all $x$
that we are considering. The monotonicity argument used for (5.13)
now shows that, if (5.14) is fulfilled, then
$$
\align
|F_N^r(x;g)|&\le C|t(-r)|^x|s(-r)|^K\frac 1{r^N}\\
&\le C|t(re^{i\epsilon})|^x|s(re^{i\epsilon})|^K\frac 1{r^N}\\
&\le C\exp\bigl[N\re
(u(re^{i\epsilon})-u(1))+y\log|t(re^{i\epsilon})|\bigr],\tag 5.21
\endalign
$$
where the last inequality is the same as in (5.19). If we take
$\epsilon=0$, we get the same right hand side as in (5.15) and hence
we obtain the same estimates, i. e. 
$$
|F_N^r(x;g)|\le C\exp\bigl[-\frac 13\min(\sqrt{\xi},\frac
14(N\rho)^{1/3})\xi \bigr].
$$
Combining this with (5.6), (5.7) and (5.18) yields
$$
|K_N(x,x)|\le CN\exp\bigl[-\frac 14\min(\sqrt{\xi},\frac
14(N\rho)^{1/3})\xi \bigr]\tag 5.22
$$
for any $\xi\ge 0$; $x$ an integer.

Consider now $\xi\in [-M_0,(\rho N)^{1/6}]$. Take $\epsilon =(\rho
N)^{-1/4}$, $\delta=\eta(\rho N)^{-1/3}\le (\rho N)^{-1/4}$, where
$\eta>0$ will be chosen below. By (5.12), we have
$$
\align
I_1'&=\frac 1{2\pi}\int_{-\epsilon}^\epsilon g((1-\delta)e^{i\theta})
\exp\bigl\{N\bigl[\frac
13\rho(1-(1-\delta)e^{i\theta})^3\\
&+
\rho (1-(1-\delta)e^{i\theta})^4v((1-\delta)e^{i\theta})\bigr]
+y\log t((1-\delta)e^{i\theta})\\
&+\omega_N\log s((1-\delta)e^{i\theta})\bigr\}d\theta.\tag 5.23
\endalign$$
We make the change of variables $\theta=\omega(\rho N)^{-1/3}$. For
$0<\eta\le (\rho N)^{1/12}$, $|\theta|\le\epsilon$, we have
$$
\frac 13\rho(1-(1-\delta)e^{i\theta})^3+
\rho (1-(1-\delta)e^{i\theta})^4v((1-\delta)e^{i\theta})=\frac
13(\eta-i\omega)^3 +R_1,\tag 5.24
$$
where $R_1\to 0$ uniformly as $N\to\infty$. Furthermore, if $\xi\in
[-M_0,(\rho N)^{1/6}]$, then
$$
y\log t((1-\delta)e^{i\theta})=(-\eta+i\omega)\xi+R_2,\tag 5.25
$$
where $R_2\to 0$ uniformly as $N\to\infty$.

Suppose $g^{(j)}(1)=0$, $j=0,\dots,\ell-1$ but $g^{(\ell)}(1)\neq 0$,
so that
$$
g((1-\delta)e^{i\theta})=\frac 1{\ell !}g^{(\ell)}(1)\rho^{-\ell/3}
(-\eta+i\omega)^\ell+\dots .\tag 5.26
$$
We now have all the estimates we need. Let $\eta=\sqrt{\xi}$ if
$\xi\ge M_0$ and $\eta=1$ if $|\xi|\le M_0$.

By (5.12) and (5.24) we obtain
$$
\re Nu((1-\delta)e^{i\theta})=\frac 13\eta^3-\eta\omega^2+R_1
$$
and hence, if $\xi\in [-M_0,(\rho N)^{1/6}]$, $\epsilon=\omega(\rho
N)^{-1/3}$ with $\omega=(\rho N)^{1/12}$,(5.19) yields, 
$$
\align
|I_1''|&\le C\exp\bigl[\frac 13\eta^3-\eta(\rho N)^{1/6}-\eta\xi+R_3\bigr]\\
&\le\frac C{N^{(\ell+1)/3}}\exp\bigl[-\frac 23|\xi|^{3/2}\bigr].\tag
5.27
\endalign
$$
Similarly, by (5.21), for $\xi\in [-M_0,(\rho N)^{1/6}]$,
$$
|I_1'|\le \frac C{N^{(\ell+1)/3}}\exp\bigl[-\frac
23|\xi|^{3/2}\bigr].\tag 5.29
$$
The dominated convergence theorem gives
$$
\align
\lim_{N\to\infty}N^{(\ell+1)/3}I_1'&=\frac{\rho^{-(\ell+1)/3}}{\ell!}
g^{(\ell)}(1)\frac 1{2\pi}\int_{-\infty}^\infty (-\eta+i\omega)^\ell
\exp\bigl[\frac i3(\omega+i\eta)^3+i\xi(\omega+i\eta)\bigr]d\omega\\
&=\frac{\rho^{-(\ell+1)/3}}{\ell!}
g^{(\ell)}(1)\Ai^{(\ell)}(\xi),\tag 5.30
\endalign
$$
uniformly for $|\xi|\le M_0$. Observe that $g_1(1)=1$,$g_2(1)=0$ but
$g_2'(1)=1$, $g_3(1)=0$ but $g_3'(1)=\rho(1-q)(\gamma q)^{-1/2}$ and
$g_4(1)=g_4'(1)=0$ but $g_4''(1)=2\rho(1-q)(\gamma
q)^{-1/2}$. Combining (5.27) and (5.29) we obtain
$$
|D_N^r(x;g)|\le\frac C{N^{(\ell+1)/3}}\exp\bigl[-\frac
23|\xi|^{3/2}\bigr],\tag 5.31
$$
for $\xi\in [-M_0,(\rho N)^{1/6}]$. The estimate (5.27) and the limit
(5.30) give
$$
\lim_{N\to\infty}N^{1/3}D_N^r(x;g_1)=\rho^{-1/3}\Ai(\xi),\tag 5.32a
$$
$$
\lim_{N\to\infty}N^{2/3}D_N^r(x;g_2)=\rho^{-2/3}\Ai'(\xi),\tag 5.32b
$$
$$
\lim_{N\to\infty}N^{2/3}D_N^r(x;g_3)=\frac{\rho^{1/3}(1-q)}{\sqrt{\gamma
q}}
\Ai'(\xi),\tag 5.32c
$$
and
$$
\lim_{N\to\infty}ND_N^r(x;g_4)=\frac{(1-q)}{\sqrt{\gamma
q}}
\Ai''(\xi),\tag 5.32d
$$
We can now use (5.22), (5.28), (5.31) and (5.32) in (5.5) and (5.6) to
prove (3.6), (3.7) and (3.8) for the Meixner kernel. The lemma is proved.

\heading 6. The equilibrium measure\endheading
\medskip
The equilibrium measure $\phi_V(t)dt$ satisfies certain variational
conditions. 

\proclaim{Proposition 6.1} Assume that $\phi\in\Cal A_s$ satisfies
\item{(i)} $\int_0^s k_V(t,\tau)\phi(\tau)d\tau\ge \lambda$ if
$\phi(t)=0$,
\item{(ii)} $\int_0^s k_V(t,\tau)\phi(\tau)d\tau\le \lambda$ if
$\phi(t)=1$, 
\item{(iii)} $\int_0^s k_V(t,\tau)\phi(\tau)d\tau = \lambda$ if
$0<\phi(t)<1$,

for some $\lambda$ (which $=F_V$). Then $\phi=\phi_V$.
\endproclaim

We will not prove this here, see [LL] for a very similar result. The
way to compute $\phi_V$ is to seek a candidate solution $\phi$ and
then verify that $\phi$ satisfies the variational conditions. 
In a region where $0<\phi(t)<1$ we can differentiate (iii) and obtain
$$
\int_0^s\frac{\phi(\tau)}{\tau -t}d\tau =-\frac 12V'(t).\tag 6.1
$$
Since $V^{\gamma,q}$ is convex the support of $\phi_V$ is a single
interval. If we consider the variational problem without the
constraint $\phi\le 1$, and this problem has a solution $\psi_0$ such
that $0\le \psi_0\le 1$, then this $\psi_0$ is the solution we are
seeking. This is the case when $\gamma\ge 1/q$, and then $[a_V,b_V]=
[a,b]$ and
$$
\int_a^b\frac{\phi(\tau)}{\tau -t}d\tau =-\frac 12V'(t),\quad
a\le t\le b.\tag 6.2
$$
We must have $\phi(b)=0$ and $\phi(a)$ bounded ($\phi(a)=0$ if
$\gamma>1/q$). 

If the solution $\psi_0(t)>1$ in some interval, e.g. $\psi_0(t)>1$ in
$[0,a_0)$ but $0<\psi_0(t)<1$ in $(a_0,b_0)$, we make an ansatz that
$\phi(t)=1$ in $[0,a]$ and $0<\phi(t)<1$ in $(a,b)$ for some $a,b$,
$[a_V,b_V]=[0,b]$. This is the situation when $\gamma<1/q$. By (6.1),
$$
\int_a^b\frac{\phi(\tau)}{\tau -t}d\tau =-\frac 12V'(t)-\int_0^a 
\frac{d\tau}{\tau -t},\tag 6.3
$$
and $\phi(a)=1$, $\phi(b)=0$. By making the substitution
$x=2(t-a)/c-1$, $y=2(\tau -a)/c-1$, $c=b-a$, in (6.2) and (6.3) we get
an equation of the form
$$
\frac 1{\pi}\int_{-1}^1\frac{v(x)}{x-y}dx=f(y),\quad -1\le x\le 1,\tag 6.4
$$
with some $f$. This equation has the general solution, [Tr],
$$
v(x)=-\frac
1{\pi\sqrt{1-x^2}}\int_{-1}^1\frac{f(y)\sqrt{1-y^2}}{y-x}dy+ 
\frac C{\pi\sqrt{1-x^2}},
$$
where $C$ is an arbitrary constant. In this way we obtain (2.19) and
(2.20).

The equation (2.21) is obtained by substituting (2.19) or (2.20) into
(2.15) (the infimum is assumed for $\tau=t$). Consider the case
$\gamma>1/q$, the other case is similar. Then, with $t=a+c(x+1)/2$,
$$
J(t)=\int_b^tJ'(s)ds=\frac c2\int_1^xJ'(a+c(y+1)/2)dy
$$
and
$$
\align
g(y)&\doteq J'(a+c(y+1)/2)=\frac c2\int_{-1}^1\frac{v(x)}{x-y}dx+
\frac 12V'(a+c(y+1)/2)\\
&=\frac c2\int_{-1}^1\log|y-x|v'(x)dx+\frac 12[\log\frac 1q-\log(y+B)+ 
\log(y+D)].
\endalign
$$
Now,
$$
v'(x)=\frac
1{2\pi}\biggl[\frac{\sqrt{D^2-1}}{x+D}-\frac{\sqrt{B^2-1}}{x+B}
\biggr]\frac 1{\sqrt{1-x^2}}
$$
and
$$
\int_{-1}^1\log|y-x|v'(x)dx=\frac 12F(y,D)-\frac 12F(y,B),
$$
where
$$
F(y,R)=\frac
1{\pi}\int_{-1}^1\frac{\sqrt{R^2-1}}{(x+R)\sqrt{1-x^2}}\log|y-x|dx.
$$
Note that
$$
\frac {d}{dy}F(y,R)=\frac{\sqrt{R^2-1}}{y+R}[\frac 1{\sqrt{y^2-1}}
+\frac 1{\sqrt{R^2-1}}].
$$
Using these formulas we see that $g(-1)=0$ and hence
$$\align
J(t)&=\frac c4\int_1^x g(y)dy=\frac c4\int_1^x(x-y)g'(y)dy\\
&=\frac
c4\int_1^x(x-y)(\frac{\sqrt{B^2-1}}{x+B}-\frac{\sqrt{D^2-1}}{x+D}
)\frac{dy}{\sqrt{y^2-1}},
\endalign
$$
which gives (2.21).

If $f(y)=(\gamma -q)(y+B)^{-1}+(1-q\gamma)(y+D)^{-1}$, then $f(y)>0$
for all $y\ge 1$ and $a_0=\inf_{1\le y\le 1/c} f(y)>0$. Thus for
$0\le\delta\le 1$, by (2.21),
$$
J(b+\delta)\ge\frac{a_0c}{8\sqrt{q\gamma}}\int_1^{1+\delta/c}
(1-\frac{2\delta}{c}-y)\frac{dy}{\sqrt{y+1}\sqrt{y-1}}\ge
c_1\delta^{3/2}, 
$$
for some constant $c_1>0$. If $\delta\ge 1$, then
$$
J(b+\delta)\ge\frac{a_0c}{8\sqrt{q\gamma}}\int_{1}^{1+1/c}
(1-\frac{2\delta}{c}-y)\frac{dy}{\sqrt{y+1}\sqrt{y-1}}
$$
which proves (2.22). A more careful computation for small $\delta$
yields (2.23).
\bigskip
\subheading{Acknowledgements} I thank C. Newman for drawing my
attention to the fact that the exponent $\chi=1/3$ occurs in many
problems. I want to express my gratitude to 
A. Dembo and P. Diaconis for telling me about the
interpretation of $G(M,N)$ in terms of randomly growing Young diagrams
and the exclusion process. Remark 2.6 was motivated by discussions
with J. Baik. Finally I thank the referee for pointing out some
mistakes in a previous version. 
This work was supported by the Swedish Natural Science Research
Council (NFR).

\bigskip
\define\kl#1{\medskip\item{[#1]    }}

\flushpar \heading REFERENCES\endheading \medskip

\kl{BR} J. Baik, E. Rains, \it Algebraic aspects of increasing
subsequences\rm, 
\newline 
math.CO/9905083

\kl{BDJ} J. Baik, P. A. Deift and K. Johansson, \it On the distribution
of the longest increasing subsequence in a random permutation\rm, 
math.CO/98101105

\kl{BDG} G. Ben Arous, A. Dembo, A. Guionnet, \it Ageing of Spherical
Spin Glasses\rm, Preprint 

\kl{BG} G. Ben Arous, A. Guionnet, \it Large Deviations for Wigner's
Law and \linebreak Voiculescu's Non-commutative Entropy\rm, Probab. Theory
Relat. Fields, {\bf 108}, (1997), 517 - 542

\kl{BO} A. Borodin, G. Olshanski, \it Statistics on Partitions, Point
Processes and the Hypergeometric Kernel\rm, math.CO/9904010

\kl{Ch} T. S. Chihara, \it An Introduction to Orthogonal
Polynomials\rm, Gordon and \linebreak Breach, New York, 1978

\kl{De} P. A. Deift, \it Orthogonal polynomials and random matrices: a 
Riemann-Hilbert approach\rm, Courant Lecture Notes in Mathematics, 3,
New York, 1999

\kl{DS} P. D. Dragnev, E. B. Saff, \it Constrained energy problems
with applications to orthogonal polynomials of a discrete variable\rm,
J. Anal. Math., {\bf 72}, (1997), 223 -  259

\kl{Fu} W. Fulton, \it Young Tableaux\rm, London Mathematical Society,
student Texts 35, Cambridge Univ. Press, 1997

\kl{Go} W. M. Y. Goh, \it Plancherel-Rotach Type Asymptotics of the
Meixner Polynomials\rm, Preprint

\kl{HP} F. Hiai, D. Petz, \it A Large Deviation Theorem for the
Empirical Eigenvalue Distribution of Random Unitary Matrices\rm,
Math. Inst. of the Hungarian Academy of Sciences, Preprint No. 17/1997

\kl{HJ} R. A. Horn, C. R. Johnson, \it Matrix Analysis\rm, Cambridge
Univ. Press, 1985

\kl{H\"o} L. H\"ormander, \it The Analysis of Linear Partial
Differential Operators I\rm, \linebreak 
Springer-Verlag, Berlin Heidelberg, 1983

\kl{Ja} A. T. James, \it Distributions of Matrix Variates and Latent
Roots Derived from Normal Samples\rm, Ann. Math. Statist., {\bf 35},
(1964), 475 - 501

\kl{JW} X.-S. Jin, R. Wang, \it Uniform Asymptotic Expansions for
Meixner Polynomials\rm, Constr. Approx., {\bf 14}, (1998), 113 - 150

\kl{JPS} W. Jockush, J. Propp, P. Shor, \it Random domino tilings and
the arctic circle theorem\rm, preprint 1995, math.CO/9801068

\kl{Jo} K. Johansson, \it On Fluctuations of Eigenvalues of Random Hermitian
Matrices\rm,  Duke Math. J., {\bf 91}, (1998), 151 - 204

\kl{Kn} D. E. Knuth, \it Permutations, Matrices and Generalized Young
Tableaux\rm, Pacific J. Math., {\bf 34}, (1970), 709 - 727

\kl{KS} J. Krug, H. Spohn, \it Kinetic Roughening of Growing
Interfaces\rm, in Solids far from Equilibrium: Growth, Morphology and
Defects , Ed. C. Godr\`eche, 479 - 582, Cambridge University Press, 1992

\kl{LL} P. D. Lax, C. D. Levermore, \it The Small Dispersion Limit of
the Korteweg-deVries Equation. I\rm, Commun. Pure and Appl. Math,
{\bf 36}, (1983), 253 - 290

\kl{Li} T. M. Ligget, \it Interacting particle Systems\rm,
Springer-Verlag, New York, 1985

\kl{NP} C. M. Newman, M. S. T. Piza, \it Divergence of Shape
Fluctuations in Two Dimensions\rm, Ann. Prob., {\bf 23}, (1995), 977 -
1005 

\kl{NSU} A. F. Nikiforov, S. K. Suslov, V. B. Uvarov, \it Classical
Orthogonal Polynomials of a Discrete Variable\rm, Springer Series in
Computational Physics, Springer-Verlag, Berlin Heidelberg, 1991

\kl{Pi} M. S. T. Piza, \it Directed Polymers in a Random Environment:
Some Results on Fluctuations\rm, J. Stat. Phys., {\bf 89}, (1997), 581
- 603 

\kl{Ro} H. Rost, \it Non-Equilibrium Behaviour of a Many Particle
Process: Density Profile and Local Equilibria\rm, Zeitschrift
f. Warsch. Verw. Gebiete, {\bf 58}, (1981), 41 - 53

\kl{ST} E. B. Saff, V. Totik, \it Logarithmic Potentials with External
Fields\rm, Grundlehren der Matematischen Wissenschaften, 316,
Springer-Verlag, Berlin, 1997

\kl{Sa} B. Sagan, \it The Symetric Group\rm, Brooks/Cole Publ. Comp., 1991

\kl{Se1} T. Sepp\"al\"ainen, \it Large Deviations for Increasing
Subsequences on the Plane\rm, Probab. Theory Relat. Fields, {\bf 112},
(1998), 221 - 244

\kl{Se2} T. Sepp\"al\"ainen, \it Coupling the totally asymmetric
simple exclusion process with a moving interface\rm, Markov
Process. Related Fields, {\bf 4}, (1998), 593- 628

\kl{Sz1} G. Szeg\"o, \it Orthogonal Polynomials\rm, American Mathematical 
Society Colloquium Publications, Vol. XXII, New York, 1939

\kl{TW1} C. A. Tracy, H. Widom, \it Level Spacing Distributions and
the Airy Kernel\rm, Commun. Math. Phys., {\bf 159}, (1994), 151 - 174 

\kl{TW2} C. A. Tracy, H. Widom, \it On Orthogonal and Symplectic
Matrix Ensembles\rm, Commun. Math. Phys., {\bf 177}, (1996), 727 - 754

\kl{TW3} C. A. Tracy, H. Widom, \it Correlation Functions, Cluster
Functions, and Spacing Distributions for Random Matrices\rm,
J. Statist. Phys., {\bf 92}, (1998), 809 - 835

\kl{Tr} F. G. Tricomi, \it Integral Equations\rm, Pure Appl. Math V,
Interscience, London, 1957

\end